\documentclass[11pt]{amsart}

\usepackage{amsmath, amsthm, amssymb, url}

\newtheorem{theo}{Theorem}[section]

\newtheorem{df}[theo]{Definition}

\newtheorem{cor}[theo]{Corollary}

\newtheorem{prop}[theo]{Proposition}

\newtheorem{claim}[theo]{Claim}

\newcommand{\R}{\mathbf{R}}
\newcommand{\C}{\mathbf{C}}
\newcommand{\Z}{\mathbf{Z}}
\newcommand{\F}{\mathbf{F}}

\newcommand{\N}{\mathbf{N}}

\newcommand{\pol}{\operatorname{pol}}
\newcommand{\Id}{\operatorname{Id}}
\newcommand{\Ad}{\operatorname{Ad}}
\newcommand{\id}{\operatorname{id}}
\newcommand{\Tr}{\operatorname{Tr}}

\newcommand{\Aut}{\operatorname{Aut}}
\newcommand{\Inn}{\operatorname{Inn}}
\newcommand{\SL}{\operatorname{SL}}
\newcommand{\alt}{\operatorname{alt}}

\newcommand{\co}{\operatorname{co}}
\newcommand{\cb}{\operatorname{cb}}
\newcommand{\Diag}{\operatorname{Diag}}
\newcommand{\op}{\operatorname{op}}
\newcommand{\dom}{\operatorname{dom}}
\newcommand{\rng}{\operatorname{rng}}
\newcommand{\bin}{\operatorname{bin}}

\begin{document}

\title[Bass-Serre rigidity results in von Neumann algebras]{Bass-Serre rigidity results in von Neumann algebras}

\begin{abstract}
We obtain new Bass-Serre type rigidity results for ${\rm II_1}$ equivalence relations and their von Neumann algebras, coming from free ergodic actions of free products of groups on the standard probability space. As an application, we show that any non-amenable factor arising as an amalgamated free product of von Neumann algebras $\mathcal{M}_1 \ast_B \mathcal{M}_2$ over an abelian von Neumann algebra $B$, is prime, i.e. cannot be written as a tensor product of diffuse factors. This gives, both in the type ${\rm II_1}$ and in the type ${\rm III}$ case, new examples of prime factors.
\end{abstract}

\author[I. Chifan]{Ionut Chifan}

\address{Math Dept \\
                 Vanderbilt University \\
                 Nashville \\
                 TN 37240 and IMAR \\
                 Bucharest \\
                 Romania}

\email{ionut.chifan@vanderbilt.edu}

\author[C. Houdayer]{Cyril Houdayer }

\address{CNRS ENS Lyon \\
UMPA UMR 5669 \\
69364 Lyon cedex 7 \\
France}

\email{cyril.houdayer@umpa.ens-lyon.fr}



\maketitle

\section{Introduction}

We prove in this paper new rigidity results for amalgamated free products (hereafter abbreviated AFP) $M = M_1 \ast_B M_2$ of semifinite von Neumann algebras over a common amenable von Neumann subalgebra. In the spirit of \cite{ipp}, these results can be viewed as von Neumann algebras analogs of ``subgroup theorems" and ``isomorphism theorems" for AFP of groups in Bass-Serre theory. Our main ``subalgebra theorem" (see Theorem \ref{subalgebra} below) shows that for any subalgebra $Q \subset M$ with no amenable direct summand (e.g. $Q$ non-amenable subfactor), the relative commutant $P = Q' \cap M$ can be conjugated by an inner automorphism of $M$ into either $M_1$ or $M_2$.

This ``subalgebra theorem" allows us to classify large classes of AFP factors as well as to prove structural properties for these algebras. For instance, we prove that any non-amenable factor $\mathcal{M} = \mathcal{M}_1 \ast_B \mathcal{M}_2$ arising as an AFP over an abelian subalgebra cannot be decomposed into a tensor product of diffuse factors, i.e. $\mathcal{M}$ is a {\em prime} factor. This gives many new examples of prime factors, of type ${\rm II_1}$ and of type ${\rm III}$. The typical ``isomorphism theorem" we derive (see Theorem \ref{bassserrevN} below) shows that if $\theta : M \simeq N^t$ is a $\ast$-isomorphism from an AFP ${\rm II_1}$ factor $M = M_1 \ast_A \cdots \ast_A M_m$ onto the amplification by some $t$ of an AFP ${\rm II_1}$ factor $N = N_1 \ast_B \cdots \ast_B N_n$, $A, B$ abelian, $1 \leq m, n \leq \infty$, with each $M_i$ and each $N_j$ containing {\em large} commuting subalgebras, then $m = n$ and $\theta(A \subset M_i)$ is unitarily conjugate to $(B \subset N_j)^t$, for all $i$, after permutation of indices.

Results of this type have been first obtained by Ozawa in his pioneering paper \cite{ozawa2004} for plain free products $M = M_1 \ast M_2$ of semi-exact ${\rm II_1}$ factors. He showed that if $Q \subset M$ is a non-amenable subfactor with $P = Q' \cap M$ hyperfinite ${\rm II_1}$ factor, then $P$ can be conjugated into $M_1$ or $M_2$ by a unitary in $M$. For the proof, he first used his C$^*$-algebraic techniques to get a ``finite dimensional" $P$-$M_i$ bimodule, and then Popa's intertwining subalgebras techniques to conclude.

In their breakthrough paper \cite{ipp}, Ioana, Peterson \& Popa showed that for $B, M_1, M_2$ arbitrary finite von Neumann algebras, any {\em rigid} subalgebra $P$ (i.e. having relative property (T)) of $M = M_1 \ast_B M_2$ can be conjugated into $M_1$ or $M_2$ by a unitary in $M$. Unlike Ozawa, this result used Popa's ``deformation/rigidity'' techniques in \cite{{popamal1}, {popa2001}} (with the relative property (T) being used for the rigidity) to first get a finite dimensional $P$-$M_i$ bimodule, then the intertwining subalgebras techniques to conclude. Peterson, using his $L^2$-derivations techniques \cite{peterson4}, obtained similar Kurosh-type results for plain free products $M = M_1 \ast M_2$, where $M_1, M_2$ are ``$L^2$-rigid" ${\rm II_1}$ factors.

In \cite{popasup}, Popa showed that
in many previous ``deformation/rigidity"
arguments 
(e.g. the W$^*$ strong rigidity of factors arising
from Bernoulli actions of property (T) groups),
the property (T) condition can be completely removed,
using instead a ``spectral gap" rigidity.
This allowed treating many
groups which do not have property (T), such as products of
an arbitrary non-amenable and
an arbitrary infinite group. Ozawa \& Popa used this ``spectral gap"
rigidity arguments \cite{ozawapopa} to prove
that the normalizer of any diffuse
amenable subalgebra of the free group factors
generates an amenable von Neumann algebra.

We combine the {\em deformation} techniques of 
\cite{{ipp}}, and the {\em spectral gap rigidity} techniques of  \cite{{popasup}, {popasolid}} in order to prove our key technical theorem. This is the above mentioned ``subalgebra
theorem", of Bass-Serre type, which is behind all the results of this paper. For $A, B \subset M$ finite von Neumann algebras, we recall from \cite{{popamal1}, {popa2001}} that the symbol $A \preceq_M B$ roughly means that a corner of $A$ can be unitarily conjugated into a corner of $B$ inside $M$.

\begin{theo}[Key intertwining result]\label{subalgebra}
For $i = 1, 2$, let $(M_i, \tau_i)$ be a finite von Neumann algebra with a common amenable von Neumann subalgebra $B \subset M_i$, such that ${\tau_1}_{|B} = {\tau_2}_{|B}$. Let $M = M_1 \ast_B M_2$ be the amalgamated free product. Let $Q \subset M$ be a von Neumann subalgebra with no amenable direct summand (e.g. $Q$ a non-amenable subfactor). Then there exists $i = 1, 2$ such that $Q' \cap M \preceq_M M_i$.
\end{theo}

We briefly recall below the concepts that we will play against each other to prove Theorem $\ref{subalgebra}$:

\begin{enumerate}
\item  The first ingredient we will use is the {\textquotedblleft malleable deformation\textquotedblright} by automorphisms $(\alpha_t, \beta)$ defined on $\tilde{M} = M \ast_B (B \overline{\otimes} L(\F_2))$, introduced in \cite{ipp}.  It represented one of the key tools that lead to the computation of the symmetry groups of AFP of weakly rigid factors. It was shown in \cite{popasup} that this deformation automatically features a certain {\textquotedblleft transversality property\textquotedblright} (see Lemma $2.1$ in \cite{popasup}) which will be of essential use in our proof.

\item The second ingredient we will use is the {\em spectral gap rigidity} principle discovered by Popa in \cite{{popasup}, {popasolid}}. We prove that for any von Neumann subalgebra $Q \subset M$ with no amenable direct summand, the action by conjugation $\Ad(\mathcal{U}(Q)) \curvearrowright \tilde{M}$ has {\textquotedblleft spectral gap\textquotedblright} relative to $M$: for any $\varepsilon > 0$, there exist $\delta > 0$ and a finite {\textquotedblleft critical\textquotedblright}  subset $F \subset \mathcal{U}(Q)$ such that for any $x \in (\tilde{M})_1$ (the unit ball of $\tilde{M}$), if $\|u x u^* - x\|_2 \leq \delta$, $\forall u \in F$, then $\|x - E_M(x)\|_2 \leq \varepsilon$.

\item A combination of $(1)$ and $(2)$ yields that for any $Q \subset M$ with no amenable direct summand, the malleable deformation $(\alpha_t)$ necessarily converges uniformly in $\| \cdot \|_2$ on the unit ball of $Q' \cap M$. Then using the intertwining techniques from \cite{ipp}, one can embed $Q' \cap M$ into $M_1$ or $M_2$ inside $M$.
\end{enumerate}

We prove in fact a more general version of Theorem $\ref{subalgebra}$  (see Theorem \ref{kurosh}), involving {\it semifinite} AFP. Indeed, given an AFP type ${\rm III}$ factor $\mathcal{M}$, these techniques
allow us to work with its {\it core} $\mathcal{M}
\rtimes_{\sigma} \R$, which is of type ${\rm II_\infty}$ \cite{{connes73}, {takesaki73}},
rather than $\mathcal{M}$ itself.
We then obtain the following theorem that generalizes many previous
results on the indecomposability of factors as tensor products, i.e. {\em primeness} (see \cite{{GaoJunge}, {Ge}, {ozawa2004}, {peterson4}, {shlya2000}}), and moreover gives new examples of prime factors, of type ${\rm II_1}$ and of type ${\rm III}$:

\begin{theo}[Primeness for AFP over abelian]\label{introresult2}
For $i = 1, 2$, let $\mathcal{M}_i$ be a von Neumann algebra. Let $B
\subset \mathcal{M}_i$ be a common abelian von Neumann subalgebra, with $B
\neq \mathcal{M}_i$, such that there exists a faithful normal
conditional expectation $E_i : \mathcal{M}_i \to B$. Denote by $\mathcal{M} =
\mathcal{M}_1 \ast_B \mathcal{M}_2$ the amalgamated free product. If
$\mathcal{M}$ is a non-amenable factor, then $\mathcal{M}$ is prime.
\end{theo}

Using some of Ueda's results on factoriality and non-amenability of plain free products and of 
amalgamated free products over a common Cartan subalgebra (see
\cite{{ueda}, {ueda3}, {ueda4}, {uedaII}}), we obtain the following corollaries:

\begin{cor}
For $i = 1, 2$, let $(\mathcal{M}_i, \varphi_i)$ be any von Neumann algebra endowed with a faithful normal state. Assume that the centralizer $\mathcal{M}_1^{\varphi_1}$ is diffuse and $\mathcal{M}_2 \neq \C$. Then the free product $(\mathcal{M}, \varphi) = (\mathcal{M}_1, \varphi_1) \ast (\mathcal{M}_2, \varphi_2)$ is a prime factor.
\end{cor}

\begin{cor}
For $i = 1, 2$, let $\mathcal{M}_i$ be a non-type ${\rm I}$ factor, and $B \subset \mathcal{M}_i$ be a common Cartan subalgebra. Then the amalgamated free product $\mathcal{M} = \mathcal{M}_1 \ast_B \mathcal{M}_2$ is a prime factor.
\end{cor}

In particular, let $\Gamma = \Gamma_1 \ast \Gamma_2$ be a free product of countable
infinite groups. Let $\sigma: \Gamma \curvearrowright (X, \mu)$ be a
free action such that the measure $\mu$ is quasi-invariant
under $\sigma$, and such that the restricted action ${\sigma}_{|\Gamma_i}$ is
 ergodic and non-transitive for $i = 1, 2$. Then the crossed product $L^\infty(X,
\mu) \rtimes \Gamma$ is a prime factor.

Theorem $\ref{subalgebra}$ allows us to obtain new W$^*$/OE Bass-Serre type rigidity results for actions of free products of groups. In order to state the main result, we first introduce the following notation. Fix integers $m, n \geq 1$. For each $i \in \{1, \dots, m\}$, and $j \in \{1, \dots, n\}$ let 
\begin{eqnarray*}
\Gamma_i & = & G_i \times H_i \\
\Lambda_j & = & G'_j \times H'_j
\end{eqnarray*}
be ICC (infinite conjugacy class) groups, such that $G_i, G'_j$ are not amenable and $H_i, H'_j$ are infinite. Denote $\Gamma = \Gamma_1 \ast \cdots \ast \Gamma_m$ and $\Lambda = \Lambda_1 \ast \cdots \ast \Lambda_n$. 

Let $\sigma : \Gamma \curvearrowright (X, \mu)$ be a free m.p. action of $\Gamma$ on the probability space $(X, \mu)$ such that $\sigma_i := \sigma_{|\Gamma_i}$ is ergodic. Write $A = L^\infty(X, \mu)$, $M_i = A \rtimes \Gamma_i$,  $M = A \rtimes \Gamma$ and $\mathcal{R}(\Gamma_i \curvearrowright X)$, $\mathcal{R}(\Gamma \curvearrowright X)$ for the associated equivalence relations.

Likewise, denote by $\rho : \Lambda \curvearrowright (Y, \nu)$ a free m.p. action of $\Lambda$ on the probability space $(Y, \nu)$ such that $\rho_j := \rho_{|\Lambda_j}$ is ergodic. Write $B = L^\infty(Y, \nu)$, $N_j = B \rtimes \Lambda_j$, $N = B \rtimes \Lambda$, and $\mathcal{R}(\Lambda_j \curvearrowright Y), \mathcal{R}(\Lambda \curvearrowright Y)$ for the associated equivalence relations. We obtain the following analogs of Theorem $7.7$ and Corollary $7.8$ of \cite{ipp}. 

\begin{theo}[W$^*$ Bass-Serre rigidity]\label{bassserrevN}
If $\theta : M \to N^t$ is a $\ast$-isomorphism, then $m = n$, $t = 1$, and after permutation of indices there exist unitaries $u_j \in N$ such that for all $j$
\begin{eqnarray*}
\Ad(u_j)\theta(M_j) & = & N_{j} \\
\Ad(u_j)\theta(A) & = & B.
\end{eqnarray*}
In particular $\mathcal{R}(\Gamma \curvearrowright X) \simeq \mathcal{R}(\Lambda \curvearrowright Y)$ and $\mathcal{R}(\Gamma_j \curvearrowright X) \simeq \mathcal{R}(\Lambda_j \curvearrowright Y)$, for any $j$.
\end{theo}

\begin{cor}[OE Bass-Serre rigidity]
If $\mathcal{R}(\Gamma \curvearrowright X) \simeq \mathcal{R}(\Lambda \curvearrowright Y)^t$, then $m = n$, $t = 1$, and after permutation of indices, we have $\mathcal{R}(\Gamma_j \curvearrowright X) \simeq \mathcal{R}(\Lambda_j \curvearrowright Y)$, for any $j$.
\end{cor}

{\bf Conventions and notations.} Throughout this paper, we write $\mathcal{M}$
for an arbitrary von Neumann algebra and $M$ for a {\it semifinite} von Neumann algebra. Usually a state is
denoted by $\varphi$ or $\psi$ and a trace is denoted by $\tau$ if
it is finite and $\Tr$ if it is semifinite. States, traces and
conditional expectations are always assumed to be faithful and
normal. We shall denote by $\mathcal{M}^n := \mathbf{M}_n(\C) \otimes \mathcal{M}$ and $\mathcal{M}^\infty := \mathbf{B}(\ell^2) \overline{\otimes} \mathcal{M}$.
Every von Neumann algebra is assumed to have separable predual.
Also, $(\mathcal{M})_1$ is the unit ball of $\mathcal{M}$
w.r.t. the operator norm.

In Section \ref{semi}, we extend some of Popa's intertwining techniques from {\it
finite} to {\it semifinite} von Neumann algebras as well as some
results from \cite{ipp}. In Section \ref{spec}, we give a generalization of Popa's {\it spectral gap} property
in the context of semifinite amalgamated free products over an
amenable von Neumann subalgebra. In Section \ref{results}, we prove
Theorem $\ref{subalgebra}$. In Section \ref{typeIII}, we prove Theorem
\ref{introresult2} and deduce several corollaries. Finally, we give
further rigidity results for finite amalgamated free products.

{\bf Acknowledgements.} The authors would like to express their warmest thanks to Stefaan Vaes for illuminating discussions on this paper. They are also very grateful to the stimulating environment at University of California, Los Angeles, where this work was done.

\section{Intertwining techniques for semifinite von Neumann algebras}\label{semi}

\subsection{Right Hilbert modules}

Let $(B, \tau)$ be a finite von Neumann algebra with a distinguished trace. Since $\tau$ is fixed, we simply denote $L^2(B, \tau)$ by $L^2(B)$. Let $H$ be a right Hilbert $B$-module, i.e. $H$ is a (separable) Hilbert space together with a normal $\ast$-representation $\pi : B^{\op} \to \mathbf{B}(H)$. For any $b \in B$, and $\xi \in H$, we shall simply denote $\pi(b^{\op}) \xi = \xi b$. By the general theory, we know that there exists an isometry $v : H \to \ell^2 \overline{\otimes} L^2(B)$ such that $v(\xi b) = v(\xi) b$, for any $\xi \in H$, $b \in B$. Since $p = vv^*$ commutes with the right $B$-action on $\ell^2 \overline{\otimes} L^2(B)$, it follows that $p \in \mathbf{B}(\ell^2) \overline{\otimes} B$. Thus, as right $B$-modules, we have $H \simeq p(\ell^2 \overline{\otimes} L^2(B))$. 

On $\mathbf{B}(\ell^2) \overline{\otimes} B$, we define the following semifinite trace $\Tr$ (which depends on $\tau$): for any $x = [x_{ij}]_{i, j} \in (\mathbf{B}(\ell^2) \overline{\otimes} B)_+$,
\begin{equation*}
\Tr \left( [x_{ij}]_{i, j} \right) = \sum_i \tau(x_{ii}).
\end{equation*}

We set $\dim(H_B) = \Tr(vv^*)$. Note that the dimension of $H$ depends on $\tau$ but does not depend on the isometry $v$. Indeed take another isometry $w : H \to \ell^2 \overline{\otimes} L^2(B)$, satisfying $w(\xi b) = w(\xi) b$, for any $\xi \in H$, $b \in B$. Note that $vw^* \in \mathbf{B}(\ell^2) \overline{\otimes} B$ and $w^*w = v^*v = 1$. Thus, we have
\begin{equation*}
\Tr(vv^*) = \Tr(v w^*w v^*) = \Tr(w v^* v w^*) = \Tr(ww^*).
\end{equation*}

Assume that $\dim(H_B) < \infty$. Then for any $\varepsilon > 0$, there exists a central projection $z \in \mathcal{Z}(B)$, with $\tau(z) \geq 1 - \varepsilon$, such that the right $B$-module $Hz$ is finitely generated, i.e. of the form $p L^2(B)^{\oplus n}$ for some projection $p \in \mathbf{M}_n(\C) \otimes B$.

\subsection{Intertwining-by-bimodules device in the semifinite setting}

In \cite{{popamal1}, {popa2001}}, Popa introduced a very powerful tool to prove the unitary conjugacy of two von Neumann subalgebras of a tracial von Neumann algebra $(M, \tau)$. If $A, B \subset (M, \tau)$ are two possibly non-unital von Neumann subalgebras, denote by $1_A, 1_B$ the units of $A$ and $B$. Note that we endow the finite von Neumann algebra $B$ with the trace $\tau(1_B \cdot 1_B) / \tau(1_B)$.

\begin{theo}[Popa, \cite{{popamal1}, {popa2001}}]\label{intertwining1}
Let $A, B \subset (M, \tau)$ be two possibly non-unital embeddings. The following are equivalent:
\begin{enumerate}
\item There exist $n \geq 1$, a possibly non-unital $\ast$-homomorphism $\psi : A \to B^n$ and a non-zero partial isometry $v \in \mathbf{M}_{1, n}(\C) \otimes 1_AM1_B$ such that $x v = v \psi(x)$, for any $x \in A$.

\item The $A$-$B$ bimodule $1_AL^2(M)1_B$ contains a non-zero $A$-$B$ subbimodule $H$ such that $\dim(H_B) < \infty$. 

\item There is no sequence of unitaries $(u_k)$ in $A$ such that
\begin{equation*}
\|E_B(a^* u_k b)\|_2 \to 0, \forall a, b \in 1_A M 1_B.
\end{equation*}
\end{enumerate}
\end{theo}
If one of the previous equivalent conditions is satisfied, we shall say that $A$ {\it embeds into} $B$ {\it inside} $M$ and denote $A \preceq_M B$.

For our purpose, we need to extend Popa's intertwining techniques for {\it semifinite} von Neumann algebras. Namely, let $(M, \Tr)$ be a von Neumann algebra endowed with a faithful normal semifinite trace. We shall simply denote by $L^2(M)$ the $M$-$M$ bimodule $L^2(M, \Tr)$, and by $\|\cdot\|_{2, \Tr}$ the $L^2$-norm associated with $\Tr$. We will use the following well-known inequality ($\|\cdot\|_\infty$ is the operator norm):
\begin{equation*}
\|x \xi y\|_{2, \Tr} \leq \|\xi\|_{2, \Tr} \|x\|_\infty \|y\|_\infty, \forall \xi \in L^2(M), \forall x, y \in M. 
\end{equation*}
We shall say that a projection $p \in M$ is $\Tr$-{\it finite} if $\Tr(p) < \infty$. Then $p$ is necessarily finite. Moreover, $pMp$ is a finite von Neumann algebra and $\tau:= \Tr(p \cdot p)/\Tr(p)$ is a faithful normal tracial state on $pMp$. Remind that for any projections $p, q \in M$, we have $p \vee q - p \sim q - p \wedge q$. Then it follows that for any $\Tr$-finite projections $p, q \in M$, $p \vee q$ is still $\Tr$-finite and $\Tr(p \vee q) = \Tr(p) + \Tr(q) - \Tr(p \wedge q)$. 

Note that if a sequence $(x_k)$ in $M$ converges to $0$ strongly, then for any non-zero $\Tr$-finite projection $q \in M$, $\left\| x_k q \right\|_{2, \Tr} \to 0$. Indeed,
\begin{eqnarray*}
x_k \to 0 \mbox{ strongly in } M & \Longleftrightarrow & x^*_kx_k \to 0 \mbox{ weakly in } M \\
& \Longrightarrow & qx^*_kx_kq \to 0 \mbox{ weakly in } qMq \\
& \Longrightarrow & \Tr(qx^*_kx_kq) \to 0 \\
& \Longrightarrow & \left\| x_k q \right\|_{2, \Tr} \to 0.
\end{eqnarray*}
Moreover, there always exists an increasing sequence of $\Tr$-finite projections $(p_k)$ in $M$ such that $p_k \to 1$ strongly.

\begin{theo}\label{intertwining}
Let $(M, \Tr)$ be a semifinite von Neumann algebra. Let $B \subset M$ be a von Neumann subalgebra such that $\Tr_{|B}$ is still semifinite. Denote by $E_B : M \to B$ the unique $\Tr$-preserving conditional expectation. Let $p \in M$ be a projection such that $\Tr(p) < \infty$. Let $A \subset pMp$ be a von Neumann subalgebra. The following conditions are equivalent:
\begin{enumerate}
\item There exists a non-zero $\Tr$-finite projection $q \in B$ such that the $A$-$qBq$ bimodule $L^2(pMq)$ contains a non-zero $A$-$qBq$ subbimodule $H$ such that $\dim(H_{qBq}) < \infty$, where $qBq$ is endowed with the trace $\Tr(q \cdot q) / \Tr(q)$.

\item There is no sequence of unitaries $(u_k)$ in $A$ such that $E_B(x^* u_k y) \to 0$ strongly, for any $x, y \in pM$.
\end{enumerate}
\end{theo}

\begin{df}
\emph{Under the assumptions of Theorem $\ref{intertwining}$, if one of the equivalent conditions is satisfied, we shall still say that $A$ {\it embeds into} $B$ {\it inside} $M$ and still denote $A \preceq_M B$.
}
\end{df}

\begin{proof}[Proof of Theorem $\ref{intertwining}$]

$(1) \Longrightarrow (2)$. Write $e = p \vee q$ which is a $\Tr$-finite projection in $M$. Thus Condition $(1)$ tells exactly that $A \preceq_{eMe} qBq$ in the sense of Theorem \ref{intertwining1}. Hence there exist $n \geq 1$, a possibly non-unital $\ast$-homomorphism $\psi : A \to (qBq)^n$ and a non-zero partial isometry $v \in M_{1, n}(\C) \otimes pMq$ such that $x v = v \psi(x)$, for any $x \in A$. Automatically, $v^*v \leq \psi(p)$.

Assume that there exists a sequence of unitaries $(u_k)$ in $A$ such that $E_B(x^* u_k y) \to 0$ strongly, for any $x, y \in pM$. In particular, we get 
\begin{equation*}
\|E_{B^n}(v^* u_k v)\|_{2, \Tr_n} \to 0.
\end{equation*}
Since $\psi(u_k)$ are unitaries in $\psi(p) B^n \psi(p)$ and $v^* u_k v  = \psi(u_k)v^*v$, we have
\begin{eqnarray*}
\|E_{B^n}(v^*v)\|_{2, \Tr_n} & = & \|E_{B^n}(v^*v) \psi(u_k)\|_{2, \Tr_n} \\
& = & \|E_{B^n}(v^*v \psi(u_k))\|_{2, \Tr_n} \\
& = & \|E_{B^n}(v^* u_k v)\|_{2, \Tr_n} \to 0.
\end{eqnarray*}
This implies that $E_{B^n}(v^*v) = 0$ and thus $v^*v = 0$. Contradiction.

$(2) \Longrightarrow (1)$. From $(2)$, there exist $\varepsilon > 0$, a finite set $F \subset L^2(B)$ and a finite set $K \subset (pM)_1$ such that
\begin{equation*}
\max_{x, y \in K, \xi \in F} \|E_B(x^* u y)\xi\| \geq \varepsilon, \forall u \in \mathcal{U}(A).
\end{equation*}
Fix an increasing sequence of $\Tr$-finite projections $(p_k)$ in $B$ such that $p_k \to 1$ strongly. Since $p_k(\xi) \to \xi$ for any $\xi \in F$, and $B \cap L^2(B)$ is dense in $L^2(B)$, we obtain that there exist $\varepsilon' > 0$ and $k_0 \in \N$ large enough, such that 
\begin{equation*}
\max_{x, y \in K} \|E_B(x^* u y)p_k\|_{2, \Tr} \geq \varepsilon', \forall u \in \mathcal{U}(A), \forall k \geq k_0.
\end{equation*}
Moreover we have
\begin{eqnarray*}
\|(1 - p_k)E_B(x^* u y)p_k\|_{2, \Tr} & \leq &  \|(1 - p_k) x^* u y p_k\|_{2, \Tr} \\
& \leq & \|(1 - p_k) x^*\|_{2, \Tr}.
\end{eqnarray*}
Since $K$ is finite and $x \in pM$ with $p$ a $\Tr$-finite projection, it follows that 
\begin{equation*}
\|(1 - p_k)E_B(x^* u y)p_k\|_{2, \Tr} \to 0,
\end{equation*}
uniformly for any $x, y \in K$, and $u \in \mathcal{U}(A)$. Thus there exist $\varepsilon'' > 0$, a $\Tr$-finite projection $q = p_{k}$ in $B$ for $k$ large enough, such that
\begin{equation*}
\max_{x, y \in K} \|qE_B(x^* u y)q\|_{2, \Tr} \geq \varepsilon'', \forall u \in \mathcal{U}(A).
\end{equation*}
The rest of the proof is now exactly the same as the one of \cite{popa2001}, because if we denote by $e = p \vee q$, we are working in the finite von Neumann algebra $e M e$.
\end{proof}

\subsection{Controlling quasi-normalizers in a semifinite AFP}

We first fix some notation. Let $P \subset Q$ be an inclusion of von Neumann algebras. We denote by
\begin{equation*}
\mathcal{N}_Q(P) := \{u \in \mathcal{U}(Q) : u P u^* = P\}
\end{equation*}
the group of all unitaries in $Q$ that {\it normalize} $P$ {\it inside} $Q$. The {\it normalizer} of $P$ {\it inside} $Q$ is the von Neumann algebra $\mathcal{N}_Q(P)''$. Since every unitary in $P' \cap Q$ normalizes $P$, we have $P' \cap Q \subset \mathcal{N}_Q(P)''$. We say that the inclusion $P \subset Q$ is {\it regular} if $\mathcal{N}_Q(P)'' = Q$. More generally, we denote by
\begin{equation*}
\mathcal{QN}_Q(P) := \{x \in Q : \exists a_1, \dots, a_n \in Q, xP \subset \sum_i P a_i \mbox{ and } Px \subset \sum_i a_i P\}
\end{equation*}
the set of all elements in $Q$ that {\it quasi-normalize} $P$ inside $Q$. Note that $\mathcal{QN}_Q(P)$ is a unital $\ast$-algebra. The {\it quasi-normalizer} of $P$ {\it inside} $Q$ is the von Neumann algebra $\mathcal{QN}_Q(P)''$. Moreover, we have $P' \cap Q \subset \mathcal{QN}_Q(P)$. We say that the inclusion $P \subset Q$ is {\it quasi-regular} if $\mathcal{QN}_Q(P)'' = Q$.

{\bf Important convention.} For $i = 1, 2$, let $M_i$ be a von Neumann algebra and $B \subset M_i$ be a common von Neumann subalgebra. Assume that there exist faithful normal conditional expectations $E_i : M \to B$. Write $(M, E) := (M_1, E_1) \ast_B (M_2, E_2)$ the amalgamated free product. We will simply denote $M = M_1 \ast_B M_2$ if no confusion is possible. We shall say that $M$ is a {\it semifinite} amalgamated free product, if there exists a semifinite faithful normal trace $\Tr$ on $M$ such that:
\begin{itemize}
\item $\Tr_{|M_i}$ and $\Tr_{|B}$ are still semifinite;
\item $\Tr \circ E = \Tr$ and $\Tr \circ E_i = \Tr$.
\end{itemize}
Whenever we refer to a trace $\Tr$ on a semifinite amalgamated free product, we always mean a trace $\Tr$ that satisfies the previous conditions.

We prove the following analog of Theorem $1.2.1$ of \cite{ipp}. Nevertheless, the proof follows the same strategy as the one of Theorem $4.6$ of \cite{houdayer3}.

\begin{theo}\label{amalga}
Let $M = M_1 \ast_B M_2$ be a semifinite amalgamated free product. Denote by $\Tr$ the semifinite trace on $M$. Let $p \in M_1$ be a projection such that $\Tr(p) < \infty$. Let $Q \subset pM_1p$ be a von Neumann subalgebra such that $Q \npreceq_{M_1} B$. Then, any $Q$-$pM_1p$ subbimodule $H$ of $L^2(pMp)$ such that $\dim(H_{pM_1p}) < \infty$, is contained in $L^2(pM_1p)$. In particular, $Q' \cap pMp$, $\mathcal{N}_{pMp}(Q)''$ and $\mathcal{QN}_{pMp}(Q)''$ are contained in $pM_1p$.
\end{theo}

\begin{proof}
Since $Q \npreceq_{M_1} B$, we know there exists a sequence of unitaries $(u_n)$ in $Q$ such that $E_B(a^* u_n b) \to 0$ strongly, for any $a, b \in pM_1$. We prove the following claim:

\begin{claim}\label{esperance}
$\forall x, y \in pMp \ominus pM_1p, \|E_{pM_1p}(x u_n y)\|_{2, \Tr} \to 0$.
\end{claim} 

\begin{proof}[Proof of Claim $\ref{esperance}$]
Let $x$ and $y$ be reduced words in $(M)_1$ with letters alternatingly from $M_1 \ominus B$ and $M_2 \ominus B$. We assume that both $x$ and $y$ contain at least a letter from $M_2 \ominus B$. We set $x = x'a$ with $a = 1$ if $x$ ends with a letter from $M_2 \ominus B$ and  $a$ equal to the last letter of $x$ otherwise. Note that $x'$ is a reduced word ending with a letter from $M_2 \ominus B$. In the same way, we set $y = by'$ with $b = 1$ if $y$ begins with a letter from $M_2 \ominus B$ and $b$ equal to the first letter of $y$ otherwise. Note that $y'$ is a reduced word beginning with a letter from $M_2 \ominus B$. Then for $z \in Q$,  we have 

\begin{equation*}
E_{M_1}(x z y) = E_{M_1}(x' E_B(a z b) y').
\end{equation*}

Apply this equality to $z = u_n$. Since $E_B(a u_n b) \to 0$ strongly, 
it follows from the equation above that $E_{M_1}(x u_n y) \to 0$ strongly. Consequently, 
\begin{equation*}
\|E_{pM_1p}(px u_n yp)\|_{2, \Tr} = \|p E_{M_1}(x u_n y) p\|_{2, \Tr} \to 0.
\end{equation*}
Using Kaplansky density theorem together with the fact that $u_n \in Q \subset p M_1 p$, we are done.
\end{proof}
The rest of the proof is now exactly the same as the one of Theorem $4.6$ in \cite{houdayer3}.
\end{proof}

\section{Spectral gap property for semifinite AFP}\label{spec}

In this section, we will freely use the language of  {\it Hilbert bimodules} over von Neumann algebras  (see \cite{{anan95}, {noncom}}). We collect here a few properties we will be using throughout. Let $M, N, P$ be any von Neumann algebras. 

For $M$-$N$ bimodules $H, K$, denote by $\pi_H$ (resp. $\pi_K$) the associated $\ast$-representation of the binormal tensor product $M \otimes_{\bin} N^{\op}$ on $H$ (resp. on $K$). We refer to \cite{EL} for the definition of $\otimes_{\bin}$. We say that $H$ is {\em weakly contained} in $K$ and denote it by $H \prec K$ if the representation $\pi_H$ is weakly contained in the representation $\pi_K$, that is if $\ker(\pi_H) \supset \ker(\pi_K)$. For a von Neumann algebra $M$, denote by $L^2(M)$ the standard representation of $M$ that gives the {\it identity} bimodule. Let $H, K$ be $M$-$N$ bimodules. The following are true:

\begin{enumerate}
\item Assume that $H \prec K$. Then, for any $N$-$P$ bimodule $L$, we have $H \overline{\otimes}_N L \prec K \overline{\otimes}_N L$, as $M$-$P$ bimodules. Exactly in the same way, for any $P$-$M$ bimodule $L$, we have $L \overline{\otimes}_M H \prec L \overline{\otimes}_M K$, as $P$-$N$ bimodules (see Lemma $1.7$ in \cite{anan95}).

\item A von Neumann algebra $B$ is amenable iff $L^2(B) \prec L^2(B) \overline{\otimes} L^2(B)$, as $B$-$B$ bimodules.
\end{enumerate}

Let $B, M, N$ be von Neumann algebras such that $B$ is amenable. Let $H$ be any $M$-$B$ bimodule and let $K$ be any $B$-$N$ bimodule. Then, as $M$-$N$ bimodules, we have $H \overline{\otimes}_B K \prec H \overline{\otimes} K$ (straightforward consequence of $(1)$ and $(2)$).

We prove the following analog of Lemma $2$ in \cite{popasolid}:

\begin{prop}\label{spectralgap}
Let $M = M_1 \ast_B M_2$ be a semifinite amalgamated free product, and denote by $\Tr$ the semifinite trace. Assume that $B$ is amenable. Let $p \in M_1$ be a projection such that $\Tr(p) < \infty$. Let $Q \subset pM_1p$ be a von Neumann subalgebra with no amenable direct summand. Then, for any free ultrafilter $\omega$ on $\N$, we have $Q' \cap (pMp)^\omega \subset (pM_1p)^\omega$.
\end{prop}

\begin{proof}
Denote $K_i = L^2(M_i) \ominus L^2(B)$, for $i = 1, 2$. Since there exists a $\Tr$-preserving normal conditional expectation $F_1 : M \to M_1$, it follows that $L^2(M_1) \ominus L^2(M_1)$ is a $M_1$-$M_1$ subbimodule of $L^2(M)$ and more precisely, we have the following isomorphism as $M_1$-$M_1$ bimodules (see \cite{{ueda}, {voiculescu92}}):
\begin{equation*}
L^2(M) \ominus L^2(M_1) \cong \bigoplus_{n \geq 1} \mathcal{H}_n,
\end{equation*}
where
\begin{equation*}
\mathcal{H}_n = L^2(M_1) \overline{\otimes}_B \mathop{\overbrace{K_2 \overline{\otimes}_B K_1 \overline{\otimes}_B \cdots \overline{\otimes}_B K_1 \overline{\otimes}_B K_2}}^{2n - 1} \overline{\otimes}_B L^2(M_1).
\end{equation*}
Let $p \in M_1$ be a non-zero $\Tr$-finite projection. Cutting down with $p$, as $pM_1p$-$pM_1p$ bimodules, we have
\begin{equation*}
L^2(pMp) \ominus L^2(pM_1p) \cong \bigoplus_{n \geq 1} p \mathcal{H}_n p.
\end{equation*}

Since $B$ is amenable, from the standard properties of composition and weak containment of correspondences recalled at the beginning of this section, it follows that as $pM_1p$-$p M_1p$ bimodules
\begin{equation*}
p \mathcal{H}_n p \prec L^2(pM_1) \overline{\otimes} \mathop{\overbrace{K_2 \overline{\otimes} K_1 \overline{\otimes} \cdots \overline{\otimes} K_1 \overline{\otimes} K_2}}^{2n - 1} \overline{\otimes} L^2(M_1p).
\end{equation*}
Consequently, we obtain
\begin{equation*}
L^2(pMp) \ominus L^2(pM_1p) \prec \bigoplus L^2(pM_1) \overline{\otimes} L^2(M_1p).
\end{equation*}
Note now that as a left $pM_1p$-module, $L^2(pM_1)$ is always a submodule of $\bigoplus L^2(pM_1p)$, and exactly the same thing for the right $pM_1p$ module $L^2(M_1p)$. Thus, we finally have
\begin{equation*}
L^2(pMp) \ominus L^2(pM_1p) \prec \bigoplus L^2(pM_1p) \overline{\otimes} L^2(pM_1p).
\end{equation*}
Since $pMp$ is a finite von Neumann algebra, the proof is then exactly the same as the one of Lemma $2$ of \cite{popasolid}.
\end{proof}

If $Q \subset pM_1p$ has no amenable direct summand, then for any $\varepsilon > 0$, there exist $\delta > 0$ and a finite subset $F \subset \mathcal{U}(Q)$ such that for any $x \in (pMp)_1$,
\begin{equation}\label{spectral}
\|ux - xu\|_{2, \Tr} < \delta, \forall u \in F \Longrightarrow \|x - E_{pM_1p}(x)\|_{2, \Tr} < \varepsilon.
\end{equation}

Remind that a ${\rm II_1}$ factor $N$ is said to be {\it full} if any central sequence is trivial, i.e. for any bounded sequence $(x_n)$ in $N$ satisfying $\|x_n y - yx_n\|_2 \to 0$ for any $y \in N$, then $\|x_n - \tau(x_n)1\|_2 \to 0$. In the case of amalgamated free products of finite von Neumann algebras, we obtain the following corollary:

\begin{cor}\label{full}
Let $(N_i, \tau_i)$ be a finite von Neumann algebra endowed with a distinguished trace, for $i = 1, 2$. Let $B \subset N_i$ be a common amenable von Neumann subalgebra such that ${\tau_1}_{|B} = {\tau_2}_{|B}$. Denote by $N = N_1 \ast_B N_2$ the amalgamated free product. If one of the $N_i$'s is a full ${\rm II_1}$ factor, then $N$ is a full ${\rm II_1}$ factor.
\end{cor}

\begin{proof}
Assume that $N_1$ is a full ${\rm II_1}$ factor. Fix $\omega$ a free ultrafilter on $\N$. We have $N_1' \cap N_1^\omega = \C$. Since $N_1$ is a non-amenable ${\rm II_1}$ factor, Proposition \ref{spectralgap} yields $N_1' \cap N^\omega \subset N_1^\omega$. Then, we obtain
\begin{eqnarray*}
N' \cap N^\omega & \subset & N_1' \cap N^\omega \\
& \subset & N_1' \cap N_1^\omega \\
& = & \C.
\end{eqnarray*}
Thus, $N$ is a full ${\rm II_1}$ factor.
\end{proof}

\section{Key intertwining theorem for semifinite AFP}\label{results}

\subsection{Notation}\label{nota}
We fix some notation that we will be using throughout this section. For $i = 1, 2$, let $N_i$ be a von Neumann algebra and let $B \subset N_i$ be a common von Neumann subalgebra such that there exist normal faithful conditional expectations $E_i : N_i \to B$. Write $N = N_1 \ast_B N_2$ the amalgamated free product. We shall always assume that $N$ is semifinite and denote by $\Tr$ its semifinite trace (see Section \ref{semi}). Write $M_i = N_i \ast_B (B \overline{\otimes} L(\Z))$ and $M = M_1 \ast_B M_2$. We still denote by $E_i : M_i \to B$ the conditional expectation, and $\Tr$ the semifinite trace on $M$. Note that $M = N \ast_B (B \overline{\otimes} L(\F_2))$. 

In $M_i$, denote by $u_i$ the generating Haar unitary of $L(\Z)$. Let $f : \mathbf{T}^1 \to ]-1, 1]$  be the  Borel function satisfying $\exp(\pi \sqrt{-1} f(z)) = z$, $\forall z \in \mathbf{T}^1$. Define $h_i = f(u_i)$ a selfadjoint element in $M_i$ such that $\exp(\pi \sqrt{-1}h_i) = u_i$. Write $u_i^t = \exp(t \pi \sqrt{-1} h_i) \in \mathcal{U}(\mathcal{M}_i)$. Following \cite{ipp}, define the deformation $(\alpha_t)$ on $M = M_1 \ast_B M_2$ by:
\begin{equation*}
\alpha_t = (\Ad u_1^t) \ast_B (\Ad u_2^t), \forall t \in \R.
\end{equation*}
Moreover define the period-$2$ automorphism $\beta$ on $M = N \ast_B (B \overline{\otimes} L(\F_2))$ by:
\begin{eqnarray*}
\beta(x) & = & x, \forall x \in N, \\
\beta(u_i) & = & u_i^*, \forall i \in \{1, 2\}. 
\end{eqnarray*}
It was proven in \cite{ipp} that the deformation $(\alpha_t)$ satisfies a certain {\it malleability type condition}:
\begin{equation*}
\alpha_t \beta = \beta \alpha_{-t}, \forall t \in \R.
\end{equation*}
Note that $\alpha_t$ and $\beta$ are $\Tr$-preserving. Hence, we shall still denote by $\alpha_t$ and $\beta$ the actions on $L^2(M)$ and note that $\beta(x) = x$, for any $x \in L^2(N)$.
Recall that the $s$-malleable deformation $(\alpha_t, \beta)$ automatically features a certain {\it tranversality property}. 
\begin{prop}[Popa, Lemma $2.1$ in \cite{popasup}]\label{transversality}
We keep the same notation as before. We have the following:
\begin{equation}\label{trans}
\|x - \alpha_{2t}(x)\|_{2, \Tr} \leq 2 \|\alpha_t(x) - E_{N}(\alpha_t(x))\|_{2, \Tr}, \; \forall x \in L^2(N), \forall t > 0.
\end{equation}
\end{prop}

\subsection{Key intertwining theorem}

All the theorems mentioned in the introduction will be consequences of the following general intertwining result:  

\begin{theo}\label{kurosh}
We keep the same notation as before. Assume that $B$ is amenable and $N = N_1 \ast_B N_2$ is a semifinite amalgamated free product, where $\Tr$ denotes the semifinite trace. Let $q \in N$ be a non-zero $\Tr$-finite projection and let $Q \subset qNq$ be a von Neumann subalgebra with no amenable direct summand. Then there exists $i \in \{1, 2\}$ such that $Q' \cap qNq \preceq_{N} N_i$.
\end{theo}

\begin{proof}
For $i \in \{1, 2\}$, write $M_i = N_i \ast_B (B \overline{\otimes} L(\Z))$, and define $M = M_1 \ast_B M_2$. Note that $M = N \ast_B (B \overline{\otimes} L(\F_2))$. We will be using the notation  introduced in subsection $\ref{nota}$. Let $q \in N$ be a non-zero projection such that $\Tr(q) < \infty$. Let $Q \subset qNq$ be a von Neumann subalgebra with no amenable direct summand. Assume that $Q' \cap qNq \npreceq_N N_i$, for all $i \in \{1, 2\}$. We shall obtain a contradiction. Denote $Q_0 = Q' \cap qNq$.

{\bf Step (1) : Using the spectral gap condition and the transversality property to find $t > 0$ and a nonzero intertwiner $v$ between $\Id$ and $\alpha_t$.}

The first step of the proof uses a well-known argument due to Popa which appeared in Theorem $4.1$ and Lemma $5.2$ in \cite{popasup}. For the sake of completeness, we will reproduce the argument from there. Let $\varepsilon = \frac{1}{8}\|q\|_{2, \Tr}$. We know that there exist $\delta > 0$, a finite subset $F \subset \mathcal{U}(Q)$, with $q \in F$, such that for every $x \in (q M q)_1$,
\begin{equation*}
\|[x, u]\|_{2, \Tr} \leq \delta, \forall u \in F \Longrightarrow \|x - E_{qNq}(x)\|_{2, \Tr} \leq \varepsilon.
\end{equation*}
Since $\alpha_t \to \Id$ pointwise $\ast$-strongly as $t \to 0$, and since for any $u \in F$
\begin{eqnarray*}
\|u - \alpha_t(u)\|^2_{2, \Tr} & = & 2 \left( \|u\|_{2, \Tr}^2 - \Re \Tr(u^*\alpha_t(u)) \right) \\
& = & 2 \left( \|u\|_{2, \Tr}^2 - \Re \Tr(q u^*\alpha_t(u)q) \right), 
\end{eqnarray*}
we may choose $t = 1/2^k$ small enough ($k \geq 1$) such that 
\begin{equation*}
\max\{\|u - \alpha_t(u)\|_{2, \Tr} : u \in F\} \leq \min \left\{ \frac{\delta}{4}, \frac{1}{48}\|q\|_{2, \Tr} \right\}.
\end{equation*}
For every $x \in Q_0$ and every $u \in F \subset Q$, writing $q = \left( q - \alpha_t(q) \right) + \alpha_t(q)$, we have
\begin{eqnarray*}
[q\alpha_t(x)q, u] & = & q\alpha_t(x)qu - uq\alpha_t(x)q \\
& = & (q - \alpha_t(q))\alpha_t(x)qu + [\alpha_t(x), u] - uq\alpha_t(x)(q - \alpha_t(q)).
\end{eqnarray*}
For every $x \in (Q_0)_1$ and every $u \in F \subset Q$, since $[u, x] = 0$, we have
\begin{eqnarray*}
\|[q\alpha_t(x)q, u]\|_{2, \Tr} & \leq & \|(q - \alpha_t(q))\alpha_t(x)qu\|_{2, \Tr} + \|[\alpha_t(x), u]\|_{2, \Tr} \\ & & + \|uq\alpha_t(x)(q - \alpha_t(q))\|_{2, \Tr} \\
& \leq & 2 \|q - \alpha_t(q)\|_{2, \Tr} + \|[\alpha_t(x), u - \alpha_t(u)]\|_{2, \Tr} \\
& \leq & 2 \|q - \alpha_t(q)\|_{2, \Tr} + 2\|u - \alpha_t(u)\|_{2, \Tr} \\
& \leq & \delta.
\end{eqnarray*}
Thus, for every $x \in (Q_0)_1$, $\|q\alpha_t(x)q - E_{qNq}(q\alpha_t(x)q)\|_{2, \Tr} \leq \varepsilon = \frac{1}{8}\|q\|_{2, \Tr}$. Now, for every $x \in Q_0$, writing $\alpha_t(q) = q + \left( \alpha_t(q) - q \right)$, we have
\begin{eqnarray*}
\alpha_t(x) & = &  q\alpha_t(x)q + (\alpha_t(q) - q) \alpha_t(x)q + q\alpha_t(x) (\alpha_t(q) - q) \\
& & + (\alpha_t(q) - q) \alpha_t(x) (\alpha_t(q) - q).
\end{eqnarray*}
Consequently, we get for every $x \in (Q_0)_1$,
\begin{eqnarray*}
\|\alpha_t(x) - E_{N}(\alpha_t(x))\|_{2, \Tr} & \leq & 6 \|q - \alpha_t(q)\|_{2, \Tr} \\
& & + \|q\alpha_t(x)q - E_{qNq}(q\alpha_t(x)q)\|_{2, \Tr} \\
& \leq & 6 \|q - \alpha_t(q)\|_{2, \Tr} + \varepsilon \\
&\leq & \frac{1}{4}\|q\|_{2, \Tr}.
\end{eqnarray*}
Using Proposition $\ref{transversality}$, we get for every $x \in (Q_0)_1$
\begin{equation*}
\|x - \alpha_s(x)\|_{2, \Tr} \leq \frac{1}{2}\|q\|_{2, \Tr},
\end{equation*}
where $s = 2t$. Thus, for every $u \in \mathcal{U}(Q_0)$, we have
\begin{eqnarray*}
\|u^*\alpha_s(u) - q\|_{2, \Tr} & = & \|u^*(\alpha_s(u) - u)\|_{2, \Tr} \\
& \leq & \|u - \alpha_s(u)\|_{2, \Tr} \\
& \leq & \frac{1}{2}\|q\|_{2, \Tr}.
\end{eqnarray*}
Denote by $\mathcal{C} = \overline{\co}^w \{u^*\alpha_s(u) : u \in \mathcal{U}(Q_0)\} \subset q M\alpha_s(q)$ the ultraweak closure of the convex hull of all $u^*\alpha_s(u)$, where $u \in \mathcal{U}(Q_0)$. Denote by $a$ the unique element in $\mathcal{C}$ of minimal $\|\cdot\|_{2, \Tr}$-norm. Since $\|a - q\|_{2, \Tr} \leq 1/2\|q\|_{2, \Tr}$, necessarily $a \neq 0$. Fix $u \in \mathcal{U}(Q_0)$. Since $u^* a \alpha_s(u) \in \mathcal{C}$ and $\|u^* a \alpha_s(u)\|_{2, \Tr} = \|a\|_{2, \Tr}$, necessarily $u^* a \alpha_s(u) = a$. Taking $v = \pol(a)$ the polar part of $a$, we have found a non-zero partial isometry $v \in qM\alpha_s(q)$ such that
\begin{equation}\label{specgap}
x v = v \alpha_s (x), \forall x \in Q_0.
\end{equation}

The rest of the proof, namely Steps $(2)$ and $(3)$, are very similar to the reasoning in Lemma $4.8$ and Theorem $6.1$ in \cite{popamsri}, Theorems $4.1$ in \cite{popamal1} and Theorem $4.3$ in \cite{ipp} (see also Theorem $5.6$ in \cite{houdayer3}). For the sake of completeness, we will nevertheless give a detailed proof.

{\bf Step (2) : Using the malleability of $(\alpha_t, \beta)$ to lift Equation $(\ref{specgap})$ till $s = 1$.}

Note that it is enough to find a non-zero partial isometry $w \in qM\alpha_{2s}(q)$ such that
\begin{equation*}
x w = w \alpha_{2s} (x), \forall x \in Q_0.
\end{equation*}
Indeed, by induction we can go till $s = 1$. Remind that $\beta(z) = z$, for every $z \in N$. Write $w = \alpha_s(\beta(v^*)v)$. Since $v \in qM\alpha_s(q)$, we see that $w \in \alpha_s\beta\alpha_s(q)M\alpha_{2s}(q)$. But $\alpha_s\beta\alpha_s = \beta$. Hence, $w \in qM\alpha_{2s}(q)$. Note that $vv^* \in Q_0' \cap qMq$. Since $Q_0 \npreceq_{N} N_i$, it follows that $Q_0 \npreceq_{N} B$. We know from Theorem $\ref{amalga}$ that $Q_0' \cap qMq \subset qNq$. In particular, $vv^* \in qNq$. Then, 
\begin{eqnarray*}
ww^* & = & \alpha_s(\beta(v^*) vv^* \beta(v)) \\
& = & \alpha_s(\beta(v^*) \beta(vv^*) \beta(v)) \\
& = & \alpha_s \beta(v^*v) \neq 0.
\end{eqnarray*}
Hence, $w$ is a non-zero partial isometry in $qM\alpha_{2s}(q)$. Moreover, for every $x \in Q_0$
\begin{eqnarray*}
w \alpha_{2s}(x) & = & \alpha_s(\beta(v^*) v \alpha_s(x)) \\
& = & \alpha_s(\beta(v^*) x v) \\
& = & \alpha_s(\beta(v^*x)v) \\
& = & \alpha_s(\beta(\alpha_s(x)v^*)v) \\
& = & \alpha_s\beta\alpha_s(x) \alpha_s(\beta(v^*)v) \\
& = & \beta(x) w \\
& = & xw.
\end{eqnarray*}

{\bf Step (3) : Using the intertwining-by-bimodules technique to conclude.}

Thus, we have found a non-zero partial isometry $v \in qM\alpha_1(q)$ such that
\begin{equation}\label{inter}
xv = v\alpha_1(x), \forall x \in Q_0.
\end{equation}
Note that $v^*v \in \alpha_1(Q_0)' \cap \alpha_1(q)M\alpha_1(q)$. Since $\alpha_1 : qMq \to \alpha_1(q)M\alpha_1(q)$ is a $\ast$-automorphism, and $Q_0 \npreceq_{N} B$,  Theorem $\ref{amalga}$ gives
\begin{eqnarray*}
\alpha_1(Q_0)' \cap \alpha_1(q)M\alpha_1(q) & = & \alpha_1\left( Q_0' \cap qMq \right) \\
& \subset & \alpha_1\left( qNq \right).
\end{eqnarray*}
Hence $v^*v \in \alpha_1(qNq)$. 

Set $A = L(\F_2)$. Denote $K_i = L^2(N_i) \ominus L^2(B)$. Denote by $P_1$ the orthogonal projection from $L^2(N)$ on $L^2(B) \oplus K_1 \oplus K_2$. Define the subspace $H_{\alt} \subset L^2(M)$ as the $L^2$-closed linear span of $B$ and the words in $N_1 \ast_B N_2 \ast_B (B \overline{\otimes} A)$ with letters alternatingly from $N_1 \ominus B$, $N_2 \ominus B$, $B \overline{\otimes} (A \ominus \C1)$ and such that two consecutive letters never come from $N_1 \ominus B$, $N_2 \ominus B$. This means that letters from $N_1 \ominus B$ and $N_2 \ominus B$ are always separated by a letter from $B \overline{\otimes} (A \ominus \C1)$.

By the definition of $\alpha_1$, it follows that $\alpha_1(L^2(N)) \subset H_{\alt}$. Denote by $P_{\alt}$ the orthogonal projection of $L^2(M)$ onto $H_{\alt}$. Since $Q_0 \npreceq_{N} N_i$, for any $i \in \{1, 2\}$, there exists a sequence of unitaries $(u_n)$ in $Q_0$ such that $E_{N_i}(a^* u_n b) \to 0$ strongly, $\forall i \in \{1, 2\}, \forall a, b \in qN$. Moreover, we have the following:

\begin{claim}\label{esperance3}
$\forall c, d \in qM\alpha_1(q)$, $\|E_{\alpha_1(N)}(c^* u_n d)\|_{2, \Tr} \to 0$.
\end{claim}

\begin{proof}[Proof of Claim $\ref{esperance3}$]
Let $c, d \in (M)_1$ (with $M = N \ast_B (B \overline{\otimes} A)$) be either in $B$ or reduced words with letters alternatingly from $N \ominus B$ and $B \overline{\otimes} (A \ominus \C1)$. Set $c = c'a$, with $a = c$ if $c \in N$, $a = 1$ if $c$ ends with a letter from $B \overline{\otimes} (A \ominus \C1)$ and $a$ equal to the last letter of $c$ otherwise. Note that $c'$ is either equal to $1$ or a reduced word ending with a letter from $B \overline{\otimes} (A \ominus \C1)$. Exactly in the same way, set $d = bd'$, with $b = d$ if $d \in N$, $b = 1$ if $d$ begins with a letter from $B \overline{\otimes} (A \ominus \C1)$ and $b$ equal to the first letter of $d$ otherwise. Note that  $d'$ is either equal to $1$ or a reduced word beginning with a letter from $B \overline{\otimes} (A \ominus \C1)$.

For $x \in N$, write $cxd = c' (axb) d'$, and note that $axb \in N$. Since $N = N_1 \ast_B N_2$ and by definition of the projection $P_{\alt}$, it is clear that 
\begin{equation*}
P_{\alt}(c' z d') = 0, \forall z \in L^2(N) \ominus (L^2(B) \oplus K_1 \oplus K_2).
\end{equation*}
By definition of the conditional expectations $E_{M_i}$, $i = 1, 2$, and $E_B$, it is easy to see that
\begin{equation*}
P_1(a u_n b)  =  E_{M_1}(a u_n b) + E_{M_2}(a u_n b) - E_B(a u_n b), \forall n.
\end{equation*}
(Note that $au_nb \in N\cap L^2(N)$.) Then $P_1(a u_n b) \to 0$ strongly. Note that by construction of $c'$ and $d'$, we have 
\begin{equation*}
P_{\alt}(c u_n d) = P_{\alt}(c' P_1(a u_n b) d').
\end{equation*}
In particular, since $\alpha_1(L^2(N)) \subset H_{\alt}$, we get
\begin{equation*}
E_{\alpha_1(N)}(c u_n d) = E_{\alpha_1(N)}(c' P_1(a u_n b) d').
\end{equation*}
Thus, $E_{\alpha_1(N)}(c u_n d) \to 0$ strongly and $\| \alpha_1(q) E_{\alpha_1(N)}(c u_n d)\alpha_1(q)\|_{2, \Tr} \to 0$. Using Kaplansky density theorem together with the fact that $u_n \in Q_0 \subset q M q$, we get the claim. 
\end{proof}

Let's come back to Equation $(\ref{inter})$. Recall that for the unitaries $(u_n)$ in $Q_0$, we have $u_n v = v \alpha_1(u_n)$. Since $\alpha_1(u_n)$ are unitaries in $\alpha_1(qNq)$, we get
\begin{eqnarray*}
\|v^*v\|_{2, \Tr} & = & \|v^*v \alpha_1(u_n)\|_{2, \Tr} \\
& = &  \|E_{\alpha_1(N)}(v^*v \alpha_1(u_n))\|_{2, \Tr} \\
& = & \|E_{\alpha_1(N)}(v^* u_n v)\|_{2, \Tr} \to 0.
\end{eqnarray*}
Hence, $v^*v = 0$. Contradiction.
\end{proof}

\section{Applications to prime factors}\label{typeIII}

\subsection{Preliminaries}
Let $\mathcal{M}$ be a von Neumann algebra. Let $\varphi$ be a state on $\mathcal{M}$. Denote by $\mathcal{M}^\varphi$ the centralizer and $M = \mathcal{M} \rtimes_{\sigma^\varphi} \R$ the {\it core} of $\mathcal{M}$, where $\sigma^\varphi$ is the modular group associated with the state $\varphi$.  
Denote by $\pi_{\sigma^\varphi} = \mathcal{M} \to M$, the representation of $\mathcal{M}$ in its core $M$, and denote by $\lambda_s$ the unitaries in $L(\R)$ implementing the action $\sigma^\varphi$. Consider the {\it dual weight} $\widehat{\varphi}$ on $M$ (see \cite{takesaki73}) which satisfies the following:
\begin{eqnarray*}
\sigma_t^{\widehat{\varphi}}(\pi_{\sigma^\varphi}(x)) & = & \pi_{\sigma^\varphi}(\sigma_t^\varphi(x)), \forall x \in \mathcal{M} \\
\sigma_t^{\widehat{\varphi}}(\lambda_s) & = & \lambda_s, \forall s \in \R.
\end{eqnarray*} 
Note that $\widehat{\varphi}$ is a semifinite weight on $M$. Write $\theta^\varphi$ the dual action of $\sigma^\varphi$ on $M$, where we identify $\R$ with its Pontryagin dual. Take now $h_\varphi$ a non-singular positive self-adjoint operator affiliated with $L(\R)$ such that $h_\varphi^{is} = \lambda_s$, for any $s \in \R$. Define $\Tr_\varphi := \widehat{\varphi}(h_\varphi^{-1} \cdot)$. We get that $\Tr_\varphi$ is a semifinite trace on $M$ and the dual action $\theta^\varphi$ {\it scales} the trace $\Tr_\varphi$:
\begin{equation*}
\Tr_\varphi \circ \theta^\varphi_s(x) = e^{-s} \Tr_\varphi(x), \forall x \in M_+, \forall s \in \R.
\end{equation*}

There is also a functorial construction of the core of the von Neumann algebra $\mathcal{M}$ which does not rely on the choice of a particular state or weight $\varphi$ (see \cite{{connes73}, {connestak}, {falcone}}). This is called the {\it non-commutaive flow of weights}. Take two states $\varphi, \psi$ on $\mathcal{M}$. It follows from \cite{falcone} and Theorem ${\rm XII.6.10}$ in \cite{takesakiII} that there exists a natural $\ast$-isomorphism
\begin{equation*}
\Pi_{\varphi, \psi} : \mathcal{M} \rtimes_{\sigma^\varphi} \R \to \mathcal{M} \rtimes_{\sigma^\psi} \R
\end{equation*}
such that $\Pi_{\varphi, \psi} \circ \theta^\varphi = \theta^\psi \circ \Pi_{\varphi, \psi}$ and $\Tr_\varphi = \Tr_\psi \circ \Pi_{\varphi, \psi}$. In the rest of this section, we will simply denote by $(M, \theta, \Tr)$ the non-commutative flow of weights, where $\theta$ is the dual action of $\R$ on the core $M$ and $\Tr$ is  the trace on $M$ such that $\Tr \circ \theta_s = e^{-s}\Tr$, for any $s \in \R$. This construction does not depend on the choice of a state on $\mathcal{M}$.

Remind that if $\mathcal{M}$ is a factor, then the dual action $\theta$ is ergodic on the center $\mathcal{Z}(M)$. We prove the following:

\begin{prop}\label{direct}
Let $\mathcal{M}$ be a non-amenable factor. Denote by $(M, \theta, \Tr)$ its non commutative flow of weights. Then for any non-zero central projection $z \in \mathcal{Z}(M)$, $Mz$ is not amenable. Moreover for any projection $p \in M$ such that $\Tr(p) < \infty$, $pMp$ is a non-amenable finite von Neumann algebra.
\end{prop}

\begin{proof}
Assume that there is a non-zero central projection $z \in \mathcal{Z}(M)$ such that $Mz$ is amenable. Thus $M\theta_t(z) = \theta_t(Mz)$ is still amenable. Define $e = \bigvee_{t \in \R} \theta_t(z)$. It is clear that $e \in \mathcal{Z}(M)$ and $\theta_t(e) = e$, for any $t \in \R$. By ergodicity of the action $\theta$ on $\mathcal{Z}(M)$, we get $e = 1$. 
Now write $z_i = \theta_{t_i}(z)$ for $i = 1, 2$, such that $Mz_i$ is amenable. Note that $z_1 \vee z_2 = z_1 + z_2 - z_1z_2 = z_1 + z_2(1 - z_1)$, so that $M(z_1 \vee z_2) = Mz_1 + Mz_2(1 - z_1)$ is still amenable. Since amenability is stable under direct limits, and since $\bigvee_{t \in \R} \theta_t(z) = 1$, it follows that $M$ is amenable. But by duality, $\mathcal{M} \overline{\otimes} \mathbf{B}(L^2(\R)) \simeq M \rtimes_\theta \R$. Consequently, $\mathcal{M}$ would be amenable. Contradiction. 

We may assume that $\mathcal{M}$ is properly infinite, so that $M$ itself is properly infinite. Let $p \in M$ be non-zero projection such that $\Tr(p) < \infty$. Denote $z = z(p)$ the central support of $p$ in $M$. Since $M$ is properly infinite, $Mz$ is still properly infinite and $Mz \simeq pMp \overline{\otimes} \mathbf{B}(\ell^2)$. Since $Mz$ is not amenable, $pMp$ is not amenable.
\end{proof}

\subsection{Main result}

We first introduce some notation. Let $(B, \tau)$ be a finite von Neumann algebra of type ${\rm I}$: for example $B = \C$, $B$ is finite dimensional or $B = L^\infty[0, 1]$. For $i = 1, 2$, let $\mathcal{M}_i$ be a von Neumann algebra endowed with a conditional expectation $E_i : \mathcal{M}_i \to B$.  We shall {\it always} assume that $B \neq \mathcal{M}_i$. Denote by $(\mathcal{M}, E) = (\mathcal{M}_1, E_1) \ast_B (\mathcal{M}_2, E_2)$ the amalgamated free product. Write 
\begin{itemize}
\item $\varphi_i = \tau \circ E_i$, $M_i = \mathcal{M}_i \rtimes_{\sigma^{\varphi_i}} \R$; 
\item $\varphi = \tau \circ E$, $M = \mathcal{M} \rtimes_{\sigma^{\varphi}} \R$.
\end{itemize}
Note that the modular groups satisfy the following equation: $\sigma_t^\varphi = \sigma_t^{\varphi_1} \ast_B \sigma_t^{\varphi_2}$, for any $t \in \R$. Set $\lambda_s$ the unitaries implementing the modular action $(\sigma_t^\varphi)$. Define the canonical conditional expectations $\widehat{E}_i : M_i \to B \overline{\otimes} L(\R)$ satisfying $\widehat{E}_i(x \lambda_s) = E_i(x) \lambda_s$, for any $x \in \mathcal{M}_i$, for any $s \in \R$. Exactly in the same way we can define $\widehat{E} : M \to B \overline{\otimes} L(\R)$. It follows from Ueda's results (see Theorem $5.1$ in \cite{ueda}) that $M$, the core of $\mathcal{M}$, is given by
\begin{equation*}
(M, \widehat{E}) = (M_1, \widehat{E}_1) \ast_{(B \overline{\otimes} L(\R))} (M_2, \widehat{E}_2).
\end{equation*}
Denote by $\Tr_{\varphi_i}$ the semifinite trace coming from the dual weight $\widehat{\varphi}_i$ on $M_i$. Exactly in the same way, denote by $\Tr_\varphi$ the semifinite trace on $M$. It is straightforward to check:
\begin{eqnarray*}
\Tr_{\varphi_i} \circ \widehat{E}_i & = &  \Tr_{\varphi_i} \\
\Tr_{\varphi} \circ \widehat{E} & = &  \Tr_{\varphi}.
\end{eqnarray*}
Moreover ${\Tr_\varphi}_{|M_i}, {\Tr_\varphi}_{|B \overline{\otimes} L(\R)}$ are still semifinite. Then $M$ is a semifinite amalgamated free product in the sense of Section \ref{semi}. We will simply denote the semifinite trace $\Tr_\varphi$ by $\Tr$ in the rest of the section.

\begin{theo}\label{prime}
Let $(B, \tau)$ be a finite von Neumann algebra of type ${\rm I}$. Assume $B \subset \mathcal{M}_i$ but $B \neq \mathcal{M}_i$. Let $E_i : \mathcal{M}_i \to B$ be a conditional expectation, for $i = 1, 2$. Denote by $\mathcal{M} = \mathcal{M}_1 \ast_B \mathcal{M}_2$ the amalgamated free product. Assume that $\mathcal{M}$ is a non-amenable factor. Then $\mathcal{M}$ is prime.
\end{theo}

\begin{proof}
We will be using the notation introduced at the beginning of this subsection. We prove the result by contradiction and we assume that $\mathcal{M}$ is not prime, i.e. $\mathcal{M} = \mathcal{P}_1 \overline{\otimes} \mathcal{P}_2$, where $\mathcal{P}_i$ is a diffuse factor (i.e. not of type ${\rm I}$). Since $\mathcal{M}$ is a non-amenable factor, we may assume that $\mathcal{P}_1$ is a non-amenable factor. Thanks to Corollary 8 of \cite{connesstormer}, we may choose a state $\psi_i$ on $\mathcal{P}_i$ such that the centralizer $\mathcal{P}_i^{\psi_i}$ is a von Neumann algebra of type ${\rm II_1}$. Denote $\psi = \psi_1 \otimes \psi_2$. Note that we can write the core $M$ in two different ways:
\begin{eqnarray*}
M & = &  (M_1, \widehat{E}_1) \ast_{(B \overline{\otimes} L(\R))} (M_2, \widehat{E}_2)\\
M & = & (\mathcal{P}_1 \overline{\otimes} \mathcal{P}_2) \rtimes_{\sigma^\psi} \R.
\end{eqnarray*}
We denote by $\Tr$ the canonical trace on $M$ scaled by the dual action $\theta$ (see the previous subsection). Denote by $P_i = \mathcal{P}_i \rtimes_{\sigma^{\psi_i}} \R$ the core of $\mathcal{P}_i$.

Fix a non-zero projection $p \in L(\R) \subset P_1$ such that $\Tr(p) < \infty$. Since $\mathcal{P}_1$ is a non-amenable factor, the finite von Neumann algebra $pP_1p$ has no amenable direct summand (see Proposition \ref{direct}). From Theorem \ref{kurosh}, we know that there exists $i = 1, 2$ such that $(pP_1p)' \cap pMp \preceq_M M_i$. Note that $\mathcal{P}_2^{\psi_2}p \subset (p P_1 p)' \cap pMp$. In particular there exists $n \geq 1$, a non-zero partial isometry $v \in \mathbf{M}_{1, n}(\C) \otimes M$, a projection $q \in M_i^n$ such that $\Tr_n(q) < \infty$ and a (unital) $\ast$-homomorphism $\rho : \mathcal{P}_2^{\psi_2}p \to q M_i^n q$ such that $x v = v \rho(x)$, for any $x \in \mathcal{P}_2^{\psi_2}p$. Denote by $Q = \langle P_1, vv^*\rangle$ the von Neumann subalgebra of $M$ generated by $P_1$ and $vv^*$. Then $Q$ is still semifinite. We have $vv^* \in Q$ and $v^*v \in \rho(\mathcal{P}_2^{\psi_2}p)' \cap qM^nq$. Since $\mathcal{P}_2^{\psi_2}p$ is a von Neumann algebra of type ${\rm II_1}$, then $\mathcal{P}_2^{\psi_2}p \npreceq_M B \overline{\otimes} L(\R)$. Consequently, Theorem \ref{amalga} yields $v^*v \in qM_i^nq$, so that we may assume $v^*v = q$. Moreover, Theorem \ref{amalga} yields $v^* Q v \subset qM_i^nq$.

With the finite projection $vv^* \in Q$, we can find a sequence of partial isometries $(u_l) \in Q$ such that $u_l^*u_l \leq vv^*$ and $z := \sum_{l} u_l u_l^* \in \mathcal{Z}(Q)$. Define $w := [ u_l v ]_l \in \mathbf{M}_{1, \infty}(\C) \overline{\otimes} M^n$. Note that, $ww^* = z$, and $w^* P_1 w \subset w^*w (M_i^n \overline{\otimes} \mathbf{B}(\ell^2)) w^*w$. But $w^*w = \Diag(v^* u^*_lu_l v)$ is not of finite trace in $M_i^n \overline{\otimes} \mathbf{B}(\ell^2)$. Note that now we are working in the semifinite amalgamated free product:
\begin{equation*}
M^\infty = M_1^\infty \ast_{(B \overline{\otimes} L(\R))^\infty} M_2^\infty.
\end{equation*}
Fix now an increasing sequence of projections $(p_k)$ in $L(\R) \subset P_1$ such that $p_0 = p$, $\Tr(p_k) < \infty$, for any $k$, and $p_k \to 1$ strongly. Denote $w_k := p_k w \in \mathbf{M}_{1, \infty}(\C) \overline{\otimes} M^n$, and note that $w_k$ is still a partial isometry (since $ww^* \in \mathcal{Z}(Q)$) and $\Tr(w_k w_k^*) < \infty$. We apply now the same strategy as before. Since $w^* \mathcal{P}_1^{\psi_1} p_k w \subset w_k^*w_k (M_i^n \overline{\otimes} \mathbf{B}(\ell^2)) w_k^*w_k$ is a subalgebra of type ${\rm II_1}$, an application of Theorem \ref{amalga} (for each $k$) yields $w^* p_k P_2 p_k w \subset w_k^*w_k(M_i^n \overline{\otimes} \mathbf{B}(\ell^2))w_k^*w_k$. Thus, we get for any $k$, 
\begin{equation*}
w^*p_k P_2 p_k w \subset w^*w(M_i^n \overline{\otimes} \mathbf{B}(\ell^2))w^*w.
\end{equation*}
Since $p_k \to 1$ strongly, we get $w^* P_2 w \subset w^*w(M_i^n \overline{\otimes} \mathbf{B}(\ell^2))w^*w$.

Consequently, we obtain $w^* P_j w \subset w^*w(M_i^n \overline{\otimes} \mathbf{B}(\ell^2))w^*w$, for any $j = 1, 2$. Since $ww^*$ commutes with $P_1$, the von Neumann algebra generated by $w^* P_1 w$ and $w^* P_2 w$ is exactly $w^* M w$, and $w^* M w \subset w^*w(M_i^{n} \overline{\otimes} \mathbf{B}(\ell^2))w^*w$. Cutting down with the projection $p_0 = p$, this implies in particular that $w^* pMp w \subset w_0^*w_0 (M_i^n  \overline{\otimes} \mathbf{B}(\ell^2)) w_0^*w_0$. Therefore, 
\begin{equation*}
w_0^*w_0 (M^n \overline{\otimes} \mathbf{B}(\ell^2)) w_0^*w_0 = w_0^*w_0 (M_i^n \overline{\otimes} \mathbf{B}(\ell^2)) w_0^*w_0.
\end{equation*}
Since $B \overline{\otimes} L(\R) \neq M_i$, by definition of the amalgamated free product $M = M_1 \ast_{(B \overline{\otimes} L(\R))} M_2$, we get a contradiction.
\end{proof}

Theorem \ref{prime} is no longer true for non-amenable factors arising as amalgamated free products over an {\it amenable} von Neumann algebra. Look at the following trivial counter-example: for $i = 1, 2$ take $N_i$ a ${\rm II_1}$ factor, write $M_i = R \overline{\otimes} N_i$, where $R$ is the hyperfinite ${\rm II_1}$ factor and $E_i = \Id \otimes \tau_i$. Then
\begin{equation*}
M := (R \overline{\otimes} N_1) \ast_R (R \overline{\otimes} N_2) = R \overline{\otimes} (N_1 \ast N_2)
\end{equation*}
is a Mc Duff ${\rm II_1}$ factor and hence not prime.

\subsection{Examples of prime factors}

We deduce now several corollaries of Theorem \ref{prime} and give new examples of prime factors. We first consider the case of {\it plain} free products. For $i = 1, 2$, let $(\mathcal{M}_i, \varphi_i)$ be any von Neumann algebra endowed with a faithful normal state. Denote by $(\mathcal{M}, \varphi) = (\mathcal{M}_1, \varphi_1) \ast (\mathcal{M}_2, \varphi_2)$ the free product. The von Neumann algebra $\mathcal{M}$ is known to be a {\it full} factor (i.e. $\Inn(\mathcal{M})$ is closed in $\Aut(\mathcal{M})$ \cite{connes74}) if one of the following conditions holds:
\begin{enumerate}
\item (Barnett \cite{barnett95}): $\exists u \in \mathcal{U}(\mathcal{M}_1^{\varphi_1}), \exists v, w \in \mathcal{U}(\mathcal{M}_2^{\varphi_2})$,
\begin{equation*}
\varphi_1(u) = \varphi_2(v) = \varphi_2(w) = \varphi_2(v^*w) = 0;
\end{equation*}

\item (Ueda \cite{{ueda3}, {ueda4}}): $\mathcal{M}_1^{\varphi_1}$ is diffuse and $\mathcal{M}_2 \neq \C$. 
\end{enumerate}
We thank Y. Ueda for pointing out to us $(2)$. Consequently, we obtain

\begin{cor}
Assume that $\mathcal{M}$ satisfies $(1)$ or $(2)$ so that $\mathcal{M}$ is a full factor. Then $\mathcal{M}$ is a prime factor.
\end{cor}

Gao \& Junge proved in \cite{GaoJunge} that any free product $(\mathcal{M}, \varphi) = \ast_{i \in I} (\mathcal{M}_i, \varphi_i)$ of amenable von Neumann algebras is solid in a general sense and hence prime. It was proven by Ricard \& Xu in \cite{RicardXu} that such a free product $\mathcal{M}$ always has the {\it complete metric approximation property}, (denoted c.m.a.p.) i.e. there exists a net $(\phi_n)$ of normal finite-rank maps on $\mathcal{M}$ such that $\limsup \|\phi_n\|_{\cb} \leq 1$ and $\phi_n \to \id_{\mathcal{M}}$ in the pointwise-ultraweak topology (see \cite{{anan95}, {CowlingHaagerup}}). Take now any countable group $\Gamma$ such that $\Lambda_{\cb}(\Gamma) > 1$, e.g. $\Gamma = \Z^2 \rtimes \SL(2, \Z)$. Then for any von Neumann algebra $\mathcal{M} \neq \C$, endowed with a faithful normal state $\varphi$, the free product $L(\Gamma) \ast (\mathcal{M}, \varphi)$ is a non-amenable factor, thus prime by Theorem \ref{prime} and which does not have the c.m.a.p. Consequently, Theorem \ref{prime} gives many examples of prime factors that do not have the c.m.a.p.

We consider now the case of {\it amalgamated} free products. Firstly, the finite case. For $i = 1, 2$, let $M_i$ be a ${\rm II_1}$ factor and $B \subset M_i$ be a common abelian von Neumann subalgebra, such that $\tau_1{_{|B}} = \tau_2{_{|B}}$. Write $M = M_1 \ast_B M_2$. We thank S. Vaes for showing us the following claim.

\begin{claim}\label{non-amenable}
The amalgamated free product $M$ is a non-amenable ${\rm II_1}$ factor.
\end{claim}

\begin{proof}[Proof of Claim $\ref{non-amenable}$]
The fact that $M$ is always a ${\rm II_1}$ factor follows from Theorem 1.1 of \cite{ipp}. We consider the following alternative:

{\bf Assume $B$ is not diffuse.} Let $p \in B$ be a non-zero minimal projection. It is straightforward to check that $pM_1p \ast_{pB} pM_2p \subset pMp$. Since $pB = \C p$, we get $pM_1p \ast pM_2p \subset pMp$. It is obvious that a free product of ${\rm II_1}$ factors is never amenable. Thus, $M$ itself is non-amenable.

{\bf Assume $B$ is diffuse.} For $n \geq 3$, since $B$ is diffuse, we may choose orthogonal projections $p_1, \dots, p_n \in B$ such that $\sum_i p_i = 1$ and $\tau(p_i) = 1/n$. Since $M_1$ and $M_2$ are both ${\rm II_1}$ factors denote by $(e_{i,j})$ (resp. $(f_{i,j})$) a system of matrix unit in $M_1$ (resp. $M_2$) such that
\begin{equation*}
e_{i, i} = f_{i, i} =  p_i, \forall i \in \{1, \dots, n\}.
\end{equation*}
Instead of writing $\{1, \dots, n\}$ for the set of indices, we will be using the notation $\Z_n := \Z/n\Z$, which is more convenient. Write
\begin{eqnarray*}
u & = & \sum_{i \in \Z_n} e_{i, i+1} \in M_1 \\
v & = & \sum_{i \in \Z_n} f_{i, i+1} \in M_2.
\end{eqnarray*}
It is straightforward to check that $u \in \mathcal{U}(M_1)$, $v \in \mathcal{U}(M_2)$ and $u^n = v^n = 1$. Moreover, we have
\begin{equation*}
u p_k u^* = v p_k v^* =  p_{k - 1}, \forall k \in \Z_n.
\end{equation*}
Hence for any $k \in \Z_n$ but $k \neq 0$ and any $j \in \Z_n$, we have $u^k p_j = p_{j - k} u^k$. Applying $E_B$, we obtain $E_B(u^k) p_j = p_{j - k} E_B(u^k) = E_B(u^k) p_{j - k}$, since $B$ is assumed to be abelian. This implies $E_B(u^k) = 0$. Likewise, we get $E_B(v^k) = 0$. It follows in particular that $u$ and $v$ are $\ast$-free in $M$ w.r.t. the trace $\tau$. Since $u$ and $v$ generate two copies of $L(\Z_n)$, we have shown that $L(\Z_n) \ast L(\Z_n) \subset M$. Since $L(\Z_n) \ast L(\Z_n)$ is non-amenable for $n \geq 3$, it follows that $M$ is a non-amenable ${\rm II_1}$ factor.
\end{proof}

\begin{cor}
For $i = 1, 2$, let $M_i$ be a ${\rm II_1}$ factor and $B \subset M_i$ be a common abelian von Neumann subalgebra. Then the amalgamated free product $M_1 \ast_B M_2$ is a prime ${\rm II_1}$ factor.
\end{cor}

More generally, for $i = 1, 2$, let now $\mathcal{M}_i$ be a non-type ${\rm I}$ factor such that $A \subset \mathcal{M}_i$ is a common {\it Cartan subalgebra}, i.e. 
\begin{itemize}
\item There exists a faithful normal conditional expectation $E_i : \mathcal{M}_i \to A$ (necessarily unique).
\item $A \subset \mathcal{M}_i$ is a MASA, i.e. $A' \cap \mathcal{M}_i = A$.
\item $A \subset \mathcal{M}_i$ is regular, i.e. $\mathcal{N}_{\mathcal{M}_i}(A)'' = \mathcal{M}_i$.
\end{itemize}
It follows from Ueda's results (see \cite{{ueda}, {uedaII}}), that under these assumptions, the amalgamated free product $\mathcal{M}_1 \ast_A \mathcal{M}_2$ is a non-amenable factor. It is even non-Mc Duff (see Theorem $8$ in \cite{ueda4}). Thus, we get

\begin{cor}
Assume that $\mathcal{M}_i$ is a non-type ${\rm I}$ factor and $A \subset \mathcal{M}_i$ is a common Cartan subalgebra. Then the amalgamated free product $\mathcal{M}_1 \ast_A \mathcal{M}_2$ is a prime factor.
\end{cor}
In particular, let $\Gamma = \Gamma_1 \ast \Gamma_2$ be a free product of countable infinite groups.  Let $\sigma : \Gamma \curvearrowright (X, \mu)$ be a free action which leaves the measure $\mu$ quasi-invariant. Assume moreover that the restriction $\sigma_{|\Gamma_i}$ is ergodic and non-transitive (see \cite{tornquist}). One can view the crossed product $\mathcal{M} := L^\infty(X, \mu) \rtimes \Gamma$ as the amalgamated free product $\mathcal{M} = \mathcal{M}_1 \ast_{L^\infty(X, \mu)} \mathcal{M}_2$, with $\mathcal{M}_i = L^\infty(X, \mu) \rtimes \Gamma_i$ and $E_i : \mathcal{M}_i \to L^\infty(X, \mu)$ is the canonical conditional expectation. Consequently, we obtain, under those assumptions, that $L^\infty(X, \mu) \rtimes \Gamma$ is a prime factor.

We point out that the assumption of ergodicity on $\sigma_{|\Gamma_i}$ cannot be removed in general. Indeed amenable factors may appear as amalgamated free products over a Cartan subalgebra. It suffices to take an amenable free ergodic action $\F_2 \curvearrowright (X, \mu)$, leaving the measure $\mu$ quasi-invariant. It follows that $\mathcal{M} = L^\infty(X, \mu) \rtimes \F_2$ is an amenable factor, hence non-prime. Nevertheless, $\mathcal{M}$ is the amalgamated free product
\begin{equation*}
\mathcal{M} = (L^\infty(X, \mu) \rtimes \Z) \ast_{L^\infty(X, \mu)} (L^\infty(X, \mu) \rtimes \Z) 
\end{equation*}

We quickly remind such a construction and refer to Section 6 of \cite{uedaII} for further details. For the free group $\F_n = \langle g_1, \dots, g_n\rangle$ on $n \geq 2$ generators, denote by $\partial\F_n$ its {\it boundary}:
\begin{equation*}
\partial\F_n = \{ (\omega_i) \in \prod_{\N^*} \{g_1, g_1^{-1}, \dots, g_n, g_n^{-1}\} : \forall i \in \N, \omega_i \neq \omega_{i + 1}^{-1} \}.
\end{equation*} 
The boundary $\partial\F_n$ is a compact space for the product topology and its topology is generated by the following clopen sets:
\begin{equation*}
\Omega(\gamma) = \{ \omega = (\omega_i) \in \partial\F_n : \omega_1 = \gamma_1, \dots, \omega_r = \gamma_r\}, 
\end{equation*}
for a reduced word $\gamma = \gamma_1 \cdots \gamma_r$. It is easy to see that the action of $\F_n$ by left multiplication on $\partial\F_n$ is continuous. By \cite{adams}, this action is known to be {\it amenable}. Consider now the probability measure $\mu$ defined on $\partial\F_n$ by:
\begin{equation*}
\mu(\Omega(\gamma)) = \frac{1}{2n}\left(\frac{1}{2n - 1}\right)^{l(\gamma) - 1}
\end{equation*}
with the word length function $l(\cdot)$. It follows from \cite{{KS}, {PS}, {RR}} that $\mu$ is quasi-invariant under $\F_n$ and moreover the action $\F_n \curvearrowright (\partial\F_n, \mu)$ is free and ergodic. It follows that the associated crossed product von Neumann algebra $L^\infty(\partial\F_n, \mu) \rtimes \F_n$ is an amenable factor. This factor is moreover of type ${\rm III_{\frac{1}{2n - 1}}}$ (see \cite{RR}). Consequently by Connes' result \cite{connes76}, it is the unique AFD factor of type ${\rm III_{\frac{1}{2n - 1}}}$, denoted by $\mathcal{R}_{\frac{1}{2n - 1}}$. Write now $\F_n = \Gamma_1 \ast \Gamma_2$ as a free product of infinite groups. Thus we have
\begin{equation*}
\mathcal{R}_{\frac{1}{2n - 1}} = (L^\infty(\partial\F_n, \mu) \rtimes \Gamma_1) \ast_{L^\infty(\partial\F_n, \mu)} (L^\infty(\partial\F_n, \mu) \rtimes \Gamma_2).
\end{equation*}

\section{Further rigidity results for finite AFP}\label{typeII}

\subsection{A Kurosh type result and consequences}
In the finite case, we can obtain some more precise results. We have the following analog of Theorem $5.1$ in \cite{ipp}:

\begin{theo}\label{kuroshII}
Let $(M_i, \tau_i)$ be a finite von Neumann algebra endowed with a distinguished trace, for $i = 1, 2$. Let $B \subset M_i$ be a common amenable von Neumann subalgebra such that ${\tau_1}_{|B} = {\tau_2}_{|B}$. Denote by $M = M_1 \ast_B M_2$ the amalgamated free product. Let $Q \subset M$ be a von Neumann subalgebra with no amenable direct summand. Denote $Q_0 = Q' \cap M$ the relative commutant. Then there exists $i = 1, 2$ such that $Q_0 \preceq_M M_i$. 

\begin{enumerate}
\item Assume that $Q_0 \npreceq_M B$. Then there exists $i = 1, 2$, $n \geq 1$ and a non-zero partial isometry $v \in M_{1, n}(\C) \otimes M$ such that $v^* Q_0 v \subset M_i^n$

\item Assume that $Q_0 \npreceq_M B$ and $M_1, M_2$ are factors. Then, there exists a unique pair of projections $q_1, q_2 \in \mathcal{Z}(Q_0' \cap M)$, satisfying $q_1 + q_2 = 1$, and unitaries $u_i \in \mathcal{U}(M)$ such that $u_i(Q_0q_i)u_i^* \subset M_i$, for $i = 1, 2$.
\end{enumerate}
\end{theo}

\begin{proof}
Theorem \ref{kurosh} yields $i = 1, 2$ such that $Q_0 \preceq_M M_i$. Thus, there exists $n \geq 1$, a projection $q \in M_i^n$, a non-zero partial isometry $v \in \mathbf{M}_{1, n}(\C) \otimes M$ and a unital $\ast$-homomorphism $\theta : Q_0 \to qM_i^nq$ such that $x v = v \theta(x)$, for any $x \in Q_0$. Assume now that $Q_0 \npreceq_M B$ and let us prove $(1)$. Since $v^*v \in \theta(Q_0)' \cap qM^nq$ and $\theta(Q_0) \npreceq_{qM^nq} B$, it follows from Theorem $1.1$ in \cite{ipp} that $v^*v \in qM_i^nq$ so that we may assume $v^*v = q$. Consequently, $v^* Q_0 v \subset M_i^n$. For $(2)$, the proof is now exactly the same as the one of Theorem $5.1$ of \cite{ipp}.
\end{proof}

Under additional assumptions, an amalgamated free product $M_1 \ast_B M_2$, with $B$ just amenable, might be prime. From a result of Hermann \& Jones (see Lemma $1$ in \cite{HJ}), if a non-inner amenable group $\Gamma$ acts on a finite von Neumann algebra $(P, \tau)$ in a trace-preserving way, then the crossed product $M = P \rtimes \Gamma$ satisfies $M' \cap M^\omega \subset P^\omega$. If we moreover assume that the action is strongly ergodic, i.e. $L(\Gamma)' \cap P^\omega = \C$, then $M' \cap M^\omega = \C$, and $M$ is a full ${\rm II_1}$ factor. Combining this observation and Corollary \ref{full}, we can obtain the following result:

\begin{theo}\label{introcor}
Let $(B, \tau)$ be any finite amenable von Neumann algebra. 
\begin{enumerate}
\item For $i = 1, 2$, let $(N_i, \tau_i)$ be finite von Neumann algebras such that ${\tau_1}_{|B} = {\tau_2}_{|B} = \tau$. Write $N = N_1 \ast_B N_2$ and assume that $N_1$ is a full ${\rm II_1}$ factor. Then $N$ is a full prime ${\rm II_1}$ factor.

\item Let $\Gamma_1, \Gamma_2$ be countable discrete groups such that $|\Gamma_1| \geq 2$ and $|\Gamma_2| \geq 3$. Denote $\Gamma = \Gamma_1 \ast \Gamma_2$, which is automatically non inner-amenable. For any strongly ergodic trace-preserving action of $\Gamma$ on $B$, $B \rtimes \Gamma$ is a full prime ${\rm II_1}$ factor. 
\end{enumerate}
\end{theo}

\begin{proof}
We note that for a non-prime ${\rm II_1}$ factor $N = N_1 \overline{\otimes} N_2$, $N_1, N_2$ are necessarily non-amenable. Then the proof is very similar to the one of Theorem \ref{kuroshII}. 
\end{proof}

\subsection{Solidity and semisolidity}
Following \cite{{ozawa2003}, {ozawa2004}}, a von Neumann algebra $M$ is said to be {\it solid} if for any diffuse von Neumann subalgebra $A \subset M$, the relative commutant $A' \cap M$ is amenable. It is said to be {\it semisolid} if for any  type ${\rm II_1}$ von Neumann subalgebra $A \subset M$, the relative commutant
 $A' \cap M$ is amenable. Ozawa proved that $L(\Gamma)$ is solid for any countable group $\Gamma$ in the class $\mathcal{S}$ (see \cite{ozawa2003}), and $L^\infty(X, \mu) \rtimes \Gamma$ is semisolid for any free ergodic m.p. action of a class $\mathcal{S}$ group $\Gamma$ on the non-atomic probability space $(X, \mu)$ (see \cite{ozawa2004}). Moreover, he showed that the following countable groups are in the class $\mathcal{S}$: word-hyperbolic groups \cite{ozawa2003}, the wreath products $\Lambda \wr \Gamma$ with $\Lambda$ amenable and $\Gamma \in \mathcal{S}$ \cite{ozawa2004} and $\Z^2 \rtimes \SL(2, \Z)$ \cite{ozawa2008}. We prove the following stability properties:
 
\begin{theo}\label{solid}
For $i = 1, 2$, let $(M_i, \tau_i)$ be a finite diffuse von Neumann algebra with a distinguished trace. 
\begin{enumerate}
\item $M_1$ and $M_2$ are solid iff the free product $M_1 \ast M_2$ is solid.
\item Take $B \subset M_i$ a common von Neumann subalgebra of type ${\rm I}$ such that ${\tau_1}_{|B} = {\tau_2}_{|B}$. Assume that $M_1$ and $M_2$ are ${\rm II_1}$ factors. Then $M_1$ and $M_2$ are semisolid iff the amalgamated free product $M_1 \ast_B M_2$ is semisolid.
\end{enumerate}
\end{theo}

\begin{proof}
Since proofs of $(1)$ and $(2)$ are similar, and since moreover $(1)$ can be deduced from Theorem $1.4$ of \cite{peterson4}, we only prove $(2)$. Since $M_1$ is a ${\rm II_1}$ factor and $B$ is of type ${\rm I}$, it follows from Theorem $1.1$ in \cite{ipp} that $M$ is a ${\rm II_1}$ factor. We prove the result by contradiction. Assume that there exists a von Neumann subalgebra $Q \subset M$ with a non-amenable direct summand such that $Q' \cap M$ is of type ${\rm II_1}$. Since $M$ is a factor, by looking at an amplification over a corner of $Q$, we may assume that $Q$ has no amenable direct summand and $Q_0 = Q' \cap M$ is still of type ${\rm II_1}$. From Theorem \ref{kuroshII}, there exist a unitary $u \in M$, a non-zero projection $q \in \mathcal{Z}(Q_0' \cap M)$, and $i = 1, 2$ such that $u (Q_0q) u^* \subset M_i$. Denote $p = uqu^*$. Theorem $1.1$ in \cite{ipp} implies that 
\begin{equation*}
u ((Q_0' \cap M)q \vee Q_0q) u^* \subset pM_ip.
\end{equation*} 
Note that $Q_0q$ is of type ${\rm II_1}$. Moreover since $Q \subset Q_0' \cap M$ and $Q$ has no amenable direct summand, it follows that $Q_0' \cap M$ has no amenable direct summand either. Thus, with $q \in \mathcal{Z}(Q_0' \cap M)$, the von Neumann algebra $(Q_0' \cap M)q$ is not amenable. This contradicts the fact that $pM_ip$ is semisolid.
\end{proof}

We cannot obtain the same statement as $(2)$ for solidity: namely, even if $M_1, M_2$ are solid and $B$ is diffuse abelian, $M = M_1 \ast_B M_2$ is not solid in general. For example take the following inclusion of free groups $\Lambda = \langle a,b^2 \rangle \subset \langle a, b \rangle = \Gamma$. Note that $[\Gamma : \Lambda] = \infty$. Look at the following generalized Bernoulli shift 
\begin{equation*}
\Gamma \curvearrowright [0, 1]^{\Gamma/\Lambda} = \Gamma \curvearrowright \prod_{g \in \Gamma/\Lambda} [0, 1]_g.
\end{equation*}
It is a free ergodic m.p. action. Write $M = L^\infty\left([0, 1]^{\Gamma/\Lambda}\right) \rtimes \Gamma$. Since $\Lambda$ acts trivially on $L^\infty([0, 1]^{e\Lambda})$ and since $\Lambda$ is not amenable, the relative commutant $L^\infty([0, 1]^{e\Lambda})' \cap M$ is not amenable.

If we want to get solidity of such an amalgamated free product, we need additional assumptions. Let $\Gamma$ be a countable group, and $\sigma : \Gamma \curvearrowright (X, \mu)$ be a free ergodic m.p. action. We shall say that this action is {\it solid} if for any diffuse von Neumann subalgebra $A \subset L^\infty(X, \mu)$, the relative commutant $A' \cap (L^\infty(X) \rtimes \Gamma)$ is amenable. We motivate this definition with the following result:

\begin{theo}
Let $\Gamma = \Gamma_1 \ast \Gamma_2$ be a free product of countable groups and consider $\Gamma \curvearrowright (X, \mu)$ a free ergodic m.p. action on the probability space. Denote by $M = L^\infty(X, \mu) \rtimes \Gamma$, $M_i = L^\infty(X, \mu) \rtimes \Gamma_i$ and note that $M = M_1 \ast_{L^\infty(X, \mu)} M_2$. Then $M$ is solid iff $M_1, M_2$ are solid and the action $\Gamma \curvearrowright (X, \mu)$ is solid.
\end{theo}

\begin{proof}
We only need to prove the {\textquotedblleft if\textquotedblright} part. We prove the result by contradiction. Since $M$ is ${\rm II_1}$ factor, there exists a von Neumann subalgebra $A \subset M$ with no amenable direct summand such that $A_0 = A' \cap M$ is diffuse. Thus, we know there exists $i = 1, 2$ such that $A_0 \preceq_M M_i$. There exist $n \geq 1$, a projection $p \in M_i^n$, a unital $\ast$-homomorphism $\psi : A_0 \to pM_i^np$ and a non-zero partial isometry $v \in \mathbf{M}_{1, n}(\C) \otimes M$ such that $x v = v \psi(x)$, for any $x \in A_0$. We may assume that $p$ equals the support projection of $E_{M_i^n}(v^*v)$.

{\bf First case: $\psi(A_0) \npreceq_{M_i^n} L^\infty(X)^n$.} The same proof as $(2)$ of Theorem \ref{solid} will lead to a contradiction, namely it will contradict the fact that $M_i$ is solid.

{\bf Second case: $\psi(A_0) \preceq_{M_i^n} L^\infty(X)^n$.} Using Remark $3.8$ in \cite{vaesbimodule}, it follows that $A_0 \preceq_M L^\infty(X)$. Then there exists $m \geq 1$, a projection $q \in L^\infty(X)^m$, a non-zero partial isometry $w \in \mathbf{M}_{1, m}(\C) \otimes M$ and a unital $\ast$-homomorphism $\theta : A_0 \to q L^\infty(X)^m q$ such that $xw = w\theta(x)$, for any $x \in A_0$. Since $\theta(A_0)$ is diffuse, by solidity of the action, $\theta(A_0)' \cap q M^m q$ is amenable. Consequently, $w^*w (\theta(A_0)' \cap q M^m q) w^*w$ is amenable and $ww^* (A_0' \cap M) ww^*$ is amenable as well. Since $A$ has no amenable direct summand, and $A \subset A_0' \cap M$, it follows that $A_0' \cap M$ has no amenable direct summand either. We get a contradiction.
\end{proof}
We refer to \cite{CI} for some applications of the notion of solid action in ergodic theory.

\subsection{W$^*$/OE Bass-Serre rigidity results}

Let $(X, \mu)$ be the standard Borel non-atomic probability space. Let $\mathcal{R}$ be a countable Borel measure-preserving equivalence relation on $(X, \mu)$. Denote by $[\mathcal{R}]$, the {\it full group} of all Borel m.p. isomorphisms $\phi : X \to X$ such that $(x, \phi(x)) \in \mathcal{R}$ for almost every $x \in X$. Denote by $[[\mathcal{R}]]$, the set of all partial Borel m.p. isomorphisms $\phi : \dom(\phi) \to \rng(\phi)$, such that $(x, \phi(x)) \in \mathcal{R}$ for almost every $x \in \dom(\phi)$. A partial Borel isomorphism $\phi \in [[\mathcal{R}]]$ is said to be {\it properly outer} if $\phi(x) \neq x$, for almost any $x \in \dom(\phi)$. Remind the following notion of {\it freeness} for equivalence relations due to Gaboriau.

\begin{df}[Gaboriau, \cite{cout}]
\emph{Let $(\mathcal{R}_k)_{k \in \N}$ be a sequence of m.p. equivalence relations on the probability space $(X, \mu)$. The sequence $(\mathcal{R}_k)$ is said to be {\it free} if for any $n \geq 1$, for any $i_1 \neq \cdots \neq i_n \in \N$, for any $\phi_j \in [[\mathcal{R}_{i_j}]]$, whenever $\phi_j$ is properly outer, the product $\phi_1 \cdots \phi_n$ is still properly outer.
}
\end{df}

In order to state the main result, we first introduce some notation. Fix integers $m, n \geq 1$. For each $i \in \{1, \dots, m\}$, and $j \in \{1, \dots, n\}$ let 
\begin{eqnarray*}
\Gamma_i & = & G_i \times H_i \\
\Lambda_j & = & G'_j \times H'_j
\end{eqnarray*}
be ICC (infinite conjugacy class) groups, such that $G_i, G'_j$ are not amenable and $H_i, H'_j$ are infinite. Note that $\Gamma_i$ and $\Lambda_j$ have a vanishing first $L^2$-Betti number (see \cite{{BV}, {PT}}). Denote $\Gamma = \Gamma_1 \ast \cdots \ast \Gamma_m$ and $\Lambda = \Lambda_1 \ast \cdots \ast \Lambda_n$. 

Let $\sigma : \Gamma \curvearrowright (X, \mu)$ be a free m.p. action of $\Gamma$ on the probability space $(X, \mu)$ such that $\sigma_i := \sigma_{|\Gamma_i}$ is ergodic. Write $A = L^\infty(X, \mu)$, $M_i = A \rtimes \Gamma_i$,  $M = A \rtimes \Gamma$, and $\mathcal{R}(\Gamma_i \curvearrowright X), \mathcal{R}(\Gamma \curvearrowright X)$ for the associated equivalence relations. 

Likewise, denote by $\rho : \Lambda \curvearrowright (Y, \nu)$ a free m.p. action of $\Lambda$ on the probability space $(Y, \nu)$ such that $\rho_j := \rho_{|\Lambda_j}$ is ergodic. Write $B = L^\infty(Y, \nu)$, $N_j = B \rtimes \Lambda_j$, $N = B \rtimes \Lambda$, and $\mathcal{R}(\Lambda_j \curvearrowright Y), \mathcal{R}(\Lambda \curvearrowright Y)$ for the associated equivalence relations. Then we have
\begin{eqnarray*}
\mathcal{R}(\Gamma \curvearrowright X) & \simeq & \mathcal{R}(\Gamma_1 \curvearrowright X) \ast \cdots \ast  \mathcal{R}(\Gamma_m \curvearrowright X) \\
\mathcal{R}(\Lambda \curvearrowright Y) & \simeq & \mathcal{R}(\Lambda_1 \curvearrowright Y) \ast \cdots \ast \mathcal{R}(\Lambda_n \curvearrowright Y).
\end{eqnarray*}

We obtain the following analogs of Theorem $7.7$ and Corollary $7.8$ of \cite{ipp}. Using our Theorem $\ref{kurosh}$, the proofs are then exactly the same. These results can be viewed as Bass-Serre type rigidity results.

\begin{theo}
If $\theta : M \to N^t$ is a $\ast$-isomorphism, then $m = n$, $t = 1$, and after permutation of indices there exist unitaries $u_j \in N$ such that for all $j$
\begin{eqnarray*}
\Ad(u_j)\theta(M_j) & = & N_{j} \\
\Ad(u_j)\theta(A) & = & B.
\end{eqnarray*}
In particular $\mathcal{R}(\Gamma \curvearrowright X) \simeq \mathcal{R}(\Lambda \curvearrowright Y)$ and $\mathcal{R}(\Gamma_j \curvearrowright X) \simeq \mathcal{R}(\Lambda_j \curvearrowright Y)$, for any $j$.
\end{theo}

\begin{cor}\label{orbitrigid}
If $\mathcal{R}(\Gamma \curvearrowright X) \simeq \mathcal{R}(\Lambda \curvearrowright Y)^t$, then $m = n$, $t = 1$, and after permutation of indices, we have $\mathcal{R}(\Gamma_j \curvearrowright X) \simeq \mathcal{R}(\Lambda_j \curvearrowright Y)$, for any $j$.
\end{cor}

Corollary \ref{orbitrigid} has been recently generalized by Alvarez \& Gaboriau \cite{AG} to all non-amenable countable groups $\Gamma_i, \Lambda_j$ with a vanishing first $L^2$-Betti number. See \cite{AG} for a precise statement.

\bibliographystyle{plain}

\begin{thebibliography}{AA}

\bibitem{adams} {\sc S. Adams}, {\it Boundary amenability for word hyperbolic groups and an application to smooth dynamics of simple groups}. Topology {\bf 33} (1994), 763--783.

\bibitem{AG} {\sc A. Alvarez \& D. Gaboriau}, {\it Free products, orbit equivalence and measure equivalence rigidity.} \url{arXiv:0806.2788}

\bibitem{anan95} {\sc C. Anantharaman-Delaroche}, {\it Amenable correspondences and approximation properties for von Neumann algebras}. Pacific J. Math. {\bf 171} (1995), 309--341.


\bibitem{barnett95} {\sc L. Barnett}, {\it Free product von Neumann algebras of type ${\rm III}$}. Proc. Amer. Math. Soc. {\bf 123} (1995), 543--553.


\bibitem{BV} {\sc B. Bekka \& A. Valette}, {\it Group cohomology, harmonic functions and the first $L^2$-Betti number}. Potential Anal. {\bf 6} (1997), 313--326.


\bibitem{CI} {\sc I. Chifan \& A. Ioana}, {\it Ergodic subequivalence relations induced by a Bernoulli action.} \url{arXiv:0802.2353}





\bibitem{noncom} {\sc A. Connes}, {\it Noncommutative geometry}. Academic Press.
San Diego, California, 1994. 

\bibitem{connes76} {\sc A. Connes}, {\it Classification of injective factors.} Ann. of Math. {\bf 104} (1976), 73--115.

\bibitem{connes74} {\sc A. Connes}, {\it Almost periodic states and factors of type ${\rm III_1}$}. J. Funct. Anal. {\bf 16} (1974), 415--445.

\bibitem{connes73} {\sc A. Connes},
{\it Une classification des facteurs de type {\rm III}.} Ann. Sci. {\'E}cole Norm. Sup. {\bf 6} (1973), 133--252.

\bibitem{connesstormer} {\sc A. Connes \& E. St\o rmer}, {\it Homogeneity of the state space of factors of type ${\rm III_1}$}. J. Funct. Anal. {\bf 28} (1978), 187--196.

\bibitem{connestak} {\sc A. Connes \& M. Takesaki},
{\it The flow of weights on factors of type {\rm III}.} T\^ohoku Math. J. {\bf 29} (1977), 473--575.

\bibitem{CowlingHaagerup} {\sc M. Cowling \& U. Haagerup},
{\it Completely bounded multipliers of the Fourier algebra of a simple Lie group of real rank one.} Invent. Math. {\bf 96} (1989), 507--549.





\bibitem{EL} {\sc E.G. Effros \& E.C. Lance}, {\it Tensor products of operator algebras.} Adv. in Math.
{\bf 25} (1977), 1--34.






\bibitem{falcone} {\sc A.J. Falcone \& M. Takesaki}, {\it Non-commutative flow of weights on a von Neumann algebra.} J. Funct. Anal. {\bf 182} (2001), 170--206.

\bibitem{cout} {\sc D. Gaboriau}, {\it Co\^ut des relations d'\'equivalence et des groupes.} Invent. Math. {\bf 139} (2000), 41--98.

\bibitem{GaoJunge} {\sc M. Gao \& M. Junge}, {\it Examples of prime von Neumann algebras.} Int. Math. Res. Notices. Vol. 2007 : article ID rnm042, 34 pages.

\bibitem{Ge} {\sc L. Ge}, {\it Applications of free entropy to finite von Neumann algebras, ${\rm II}$.} Ann. of Math. {\bf 147} (1998), 143--157.

\bibitem{HJ} {\sc R.H. Hermann \& V.F.R. Jones}, {\it Central sequences in crossed products.} Cont. Math. {\bf 62} (1987), 539--544.

\bibitem{houdayer3} {\sc C. Houdayer}, {\it Construction of type ${\rm II_1}$ factors with prescribed countable fundamental group.} J. reine angew Math. {\bf 634} (2009), 169-207.



\bibitem{ipp} {\sc A. Ioana, J. Peterson \& S. Popa}, {\it Amalgamated free products of $w$-rigid factors and calculation of their symmetry groups.} Acta Math. {\bf 200} (2008), 85--153. 





\bibitem{KS} {\sc G. Kuhn \& T. Steger}, {\it More irreducible boundary representations of free
groups.} Duke Math. J. {\bf 82} (1996), 381--436.

\bibitem{ozawa2003} {\sc N. Ozawa}, {\it Solid von Neumann algebras.} Acta Math. {\bf 192} (2004), 111--117. 

\bibitem{ozawa2004} {\sc N. Ozawa} {\it A Kurosh-type theorem for type ${\rm II_1}$ factors.} Int. Math. Res. Notices. Vol. 2006 : article ID 97560, 21 pages.

\bibitem{ozawa2008} {\sc N. Ozawa}, {\it An example of a solid von Neumann algebra.} Hokkaido Math. J. {\bf 38} (2009), 557--561. 

\bibitem{ozawapopa} {\sc N. Ozawa \& S. Popa}, {\it On a class of $\rm{II}_1$ factors with at most one Cartan subalgebra.} Ann. of Math., to appear. \url{arXiv:0706.3623} 

\bibitem{PS} {\sc C. Pensavalle \& T. Steger}, {\it Tensor products with anisotropic principal series
representations of free groups.} Pacific J. Math. {\bf 173} (1996), 181--202.

\bibitem{peterson4} {\sc J. Peterson}, {\it $L^2$-rigidity in von Neumann algebras.} Invent. Math. {\bf 175} (2009), 417--433.

\bibitem{PT} {\sc J. Peterson \& A. Thom}, {\it Group cocycles and the ring of affiliated operators.} \url{arXiv:0708.4327}

\bibitem{popasup} {\sc S. Popa}, {\it On the superrigidity of malleable actions with spectral gap.}  J. Amer. Math. Soc. {\bf 21} (2008), 981--1000. 

\bibitem{popasolid} {\sc S. Popa}, {\it On Ozawa's property for free group factors.} Int. Math. Res. Notices.  Vol. 2007 : article ID rnm036, 10 pages.


\bibitem{popamal1} {\sc S. Popa}, {\it Strong rigidity of ${\rm II_1}$ factors arising from malleable actions of w-rigid groups ${\rm I}$.} Invent. Math. {\bf 165} (2006), 369-408.


\bibitem{popa2001} {\sc S. Popa}, {\it On a class of type ${\rm II_1}$ factors with Betti numbers invariants.} Ann. of Math. {\bf 163} (2006), 809--899.

\bibitem{popamsri} {\sc S. Popa}, {\it Some rigidity results for non-commutative Bernoulli Shifts.} J. Funct. Anal. {\bf 230} (2006), 273--328.




\bibitem{RR} {\sc J. Ramagge \& G. Robertson}, {\it Factors from trees.} Proc. Amer. Math. Soc.
{\bf 125} (1997), 2051--2055.


\bibitem{RicardXu} {\sc \'E. Ricard \& Q. Xu}, {\it Khintchine type inequalities for reduced free products and applications.} J. reine angew. Math. {\bf 599} (2006), 27--59.




\bibitem{shlya2000} {\sc D. Shlyakhtenko}, {\it Prime type ${\rm III}$ factors.}  Proc. Nat. Acad. Sci., {\bf 97} (2000), 12439--12441.





\bibitem{takesakiII} {\sc M. Takesaki}, { \it Theory of Operator Algebras ${\rm II}$.} EMS {\bf 125}. Springer-Verlag, Berlin, Heidelberg, New-York, 2000.

\bibitem{takesaki73} {\sc M. Takesaki}, {\it Duality for crossed products and structure of von Neumann algebras of type ${\rm III}$.} Acta Math. {\bf 131} (1973), 249--310.

\bibitem{tornquist} {\sc A. T\" ornquist}, {\it Orbit equivalence and actions of $\F_n$.} J. Symbolic Logic {\bf 71} (2006), 265--282.

\bibitem{ueda} {\sc Y. Ueda}, {\it Amalgamated free products over Cartan subalgebra.} Pacific J. Math. {\bf 191} (1999), 359--392.

\bibitem{ueda3} {\sc Y. Ueda}, {\it Remarks on free products with respect to non-tracial states.} Math. Scand. {\bf 88} (2001), 111--125.

\bibitem{ueda4} {\sc Y. Ueda}, {\it Fullness, Connes' $\chi$-groups, and ultra-products of amalgamated free products over Cartan subalgebras.} Trans. Amer. Math. Soc. {\bf 355} (2003), 349--371.

\bibitem{uedaII} {\sc Y. Ueda}, {\it Amalgamated free products over Cartan subalgebra, ${\rm II}$. Supplementary results and examples.} Advanced Studies in Pure Mathematics, Operator Algebras and Applications. {\bf 38} (2004),  239--265.


\bibitem{vaesbern} {\sc S. Vaes}, {\it Rigidity results for Bernoulli actions and their von Neumann algebras (after S. Popa).} S{\'e}minaire Bourbaki, expos\'e 961. Ast\'erisque {\bf 311} (2007), 237-294.

\bibitem{vaesbimodule} {\sc S. Vaes}, {\it Explicit computations of all finite index bimodules for a family of ${\rm II_1}$ factors.} Ann. Sci. \'Ecole Norm. Sup. {\bf 41} (2008), 743--788.




\bibitem{voiculescu92} {\sc D.-V. Voiculescu, K.J. Dykema \& A. Nica}, {\it Free
  random variables.} CRM Monograph Series {\bf 1}.
American Mathematical Society, Providence, RI, $1992$. 




\end{thebibliography}

\end{document}